\magnification \magstep0
\input amssym.tex\input amstex
\loadmsam
\loadmsbm
\loadbold

\documentstyle{amsppt}

\magnification = \magstep0

\def\rtzero{\rho}

\def\azero{a}
\def\uhat{{\hat u}}
\def\vhat{{\hat v}}
\def\zhat{{\hat z}}

\def\div{\hbox{div}}
\def\Pic{\hbox{Pic}^0(C)}
\def\Jac{\hbox{Jac}(C)}
\def\Caff{C_{\hbox{\sevenrm aff}}}
\def\nl{\cr&\qquad }

\def\bolds{{\bold s}}
\def\boldhats{{\hat \bold s}}
\def\boldS{{\bold S}}
\def\boldhatS{{\hat \bold S}}
\def\bolda{{\bold a}}
\def\boldazero{{\bold a}_0}

\def\boldg{{\bold g}}
\def\bolde{{\bold e}}
\def\boldG{{\bold G}}
\def\boldE{{\bold E}}
\def\boldzero{{\bold 0}}
\def\boldx{{\bold x}}
\def\boldx{{\bold x}}

\def\boldb{{\bold b}}
\def\boldc{{\bold c}}
\def\boldd{{\bold d}}

\def\bolda{{\bold a}}

\def\Z{\Bbb Z}

\def\P{\Bbb P}
\def\A{\Bbb A}
\def\H{\Cal H}

\def\abs#1{\vert{#1}\vert}

\input V3,95
\newarrow{Corresponds}<--->
\newarrow{Twoway}<--->
\newarrow{Backwards}<----

\def\cubesize{1.3em}
\def\squaresize{1.9em}

\def\labeledsquare#1#2#3#4#5#6#7#8{
\diagramstyle[size=\squaresize]
\diagram
{#1} & \rCorresponds^{#5} & {#2} \\
\dCorresponds^{#6} & & \dCorresponds^{#7} \\
{#3} & \rCorresponds^{#8} & {#4} \\
\enddiagram
}

\def\cube#1#2#3#4#5#6#7#8{
\diagramstyle[size=\cubesize]
\diagram
{#5}     &             & \rTwoway &          & {#6}        &           &          \\
         & \rdTwoway   &          &          & \dBackwards & \rdTwoway &          \\
\dTwoway &             & {#1}     & \rTwoway & \HonV       &           & {#2}     \\
         &             & \dTwoway &          & \dTo        &           &          \\
{#7}     & \rBackwards & \VonH    & \rTo     & {#8}        &           & \dTwoway \\
         & \rdTwoway   &          &          & & \rdTwoway &           &          \\
         &             & {#3}     &          & \rTwoway    &           & {#4}     \\
\enddiagram}

\def\cubecorner#1#2#3#4{
\diagramstyle[size=\cubesize]
\diagram
{#4} &           &          &          & &      \\
     & \rdTwoway &          &          & &      \\
     &           & {#1}     & \rTwoway & & {#2} \\
     &           & \dTwoway &          & &      \\
     &           &          &          & &      \\
     &           & {#3}     &          & &      \\
\enddiagram}

\def\cubeminuscorner#1#2#3#4#5#6#7{
\diagramstyle[size=\cubesize]
\diagram
{#5} & & \rTwoway & & {#6} & & \\
& \rdTwoway & & & & \rdTwoway & \\
\dTwoway & & {#1} & \rTwoway & & & {#2} \\
& & \dTwoway & & & & \\
{#7} & & & & & & \dTwoway \\
& \rdTwoway & & & & & \\
& & {#3} & & \rTwoway & & {#4} \\
\enddiagram}

\def\zz{00}
\def\zo{01}
\def\oz{10}
\def\oo{11}

\def\zzz{000}
\def\zzo{001}
\def\zoz{010}
\def\zoo{011}
\def\ozz{100}
\def\ozo{101}
\def\ooz{110}
\def\ooo{111}

\def\I{{\Cal I}}

\def\OR{\hbox{ OR }}
\def\AND{\hbox{ AND }}

\def\rtone{\rho^{(1)}}

\def\boldaone{{\bold a}^{(1)}}

\def\rttwo{\rho^{(2)}}

\def\rti{\rho^{(i)}}

\def\boldz{{\bold z}}
\def\JacC{\hbox{Jac}(C)}
\def\Cl{\hbox{Cl}}

\title
Mult-Projective Models for Jacobian Varieties of Genus Two Curves
\endtitle

\author
Mark Heiligman
\endauthor

\address
Intelligence Advanced Research Projects Activity,
College Park, Maryland
\endaddress

\email
mark.heiligman\@iarpa.gov
\endemail

\date
January 1, 2015
\enddate


\abstract 
A set of multi-homogeneous equations for the Jacobian of a genus two
curve is given.  The approach used is to write down affine equations
for the Jacobian minus various tranlations of the Theta-divisor by
[2]-division points, and then to write down affine glue equations for
the overlaps. Taking multi-projective completions for all of these
then yields a complete multi-projective two-dimensional variety whose
points are in one-to-one correspondence with degree zero divisor
classes on the curve (i.e.  the Picard group). This multi-projective
variety then becomes a complete projective two-dimensional variety
under the Segre imbedding.
\endabstract 

\rightheadtext{Genus 2 Jacobians}

\thanks
{\it Disclaimer}. 
All statements of fact, opinion, or analysis expressed in this paper
are solely those of the author and do not necessarily reflect the
official positions or views of the Office of the Director of National
Intelligence (ODNI), the Intelligence Advanced Research Projects
Activity (IARPA), or any other government agency.  Nothing in the
content should be construed as asserting or implying U.S. Government
authentication of information or ODNI endorsement of the author's
views.
\endthanks

\endtopmatter

\document

\head
Introduction and Motivation
\endhead

The Jacobian variety associated with an algebraic curve is well known
and widely studied. For such an important mathematical object, it is
remarkable how difficult it is to write out defining equations. Of
course, in genus one, which are elliptic curves, there is no
problem. However, for anything beyond the genus one case, the
situation is vastly more difficult. Even for genus two, where the
Jacobian is an abelian surface, a set of defining equations is
daunting. 

From a fairly elementary viewpoint, the group of divisor classes of
degree 0 on $C$ is the Picard group $\Pic$. This is defined as an
abstract group and {\it a priori} does not have the structure of an
algebraic variety. The realization of $\Pic$ as a projective algebraic
variety is simply $\Jac$.

Few authors even try to write equations for abelian varieties (see
site{Mi1}, \cite{Mi2}, and \cite{Mu1}, for example) and even when
they do, (see \cite{Mu2}), the approaches are far from concrete.
Previous works by Flynn (see \cite{Fl1}, \cite{Fl2}, and \cite{CF}),
Grant (see \cite{Gr}), and Wamellen (\cite{Wa}) have yielded sets of
such defining equations, but they tend to be long and
complicated. Their approach has been to write down a basis for the
linear system of a very ample divisor on the Jacobian and then take
all the relations between functions in this linear system. Another
approach by Anderson (see \cite{An}) also yields a complicated set of
equations based on invariant theory. 

There are many reasons for the difficulty of writing out such a set of
equations. As has been noted by Mumford and Cantor (see \cite{Ca} and
\cite{Mu3}), building on an approach going back to Jacobi, writing out
a set of defining equations for an affine variety for a large piece of
the Jacobian of a hyperelliptic curve $C$ is relatively
straightforward. In general, if $C$ is an algebraic curve of genus
$g$, then $\Jac$ is a projective variety (and therefore complete) of
dimension $g$.

However, $\Jac-\Theta$ is an affine variety, where $\Theta$ (the
theta divisor) is a subvariety of codimension 1, and in the case of
$g=2$, the variety $\Theta$ is simply a copy of the curve $C$ imbedded
into $\Jac$. It is relatively straightforward to show
that $\Jac-\Theta$ can be written as a system of $g$ equations in $2g$
variables. The problem is that the projective completion of
$\Jac-\Theta$ does not yield $\Jac$. In fact, for $g>1$, $\Jac$ is not
a complete intersection, so there will always be an overdetermined
system of defining equations for $\Jac$.  (One way of seeing this
problem is that the middle dimension of the cohomology of $\Jac$ is
non-zero, while for a complete intersection, only the highest
dimension cohomology is non-zero.)

The goal of this note is to show how to write out a set of equations
for $\Jac$ that are not particularly complicated. The idea is to
realize $\Jac$ as a multi-projective variety.  Then the Segre imbedding
of multi-projective space into ordinary projective space will give
$\Jac$ as a projective variety.

This approach is implicit in Mumford's work, although it is not
explicitly worked out there. The idea is that multiple copies of the
affine variety $\Jac-\Theta$ can be glued together to yield a
projective variety whose points are in one-to-one correspondence with
the elements of $\Pic$. These glue equations are also fairly simple
and simply represent the results of adding points of order 2 in
$\Pic$. More specifically, what really gets glued together are
multiple copies of the projectivization of $\Jac-\Theta$, and what
needs to be carefully checked is that there are no ``extraneous''
points on this glued together variety, i.e. that every geometric point
corresponds to an actual element of $\Pic$.

One drawback of this approach is that the projective space in which
$\Jac$ is imbedded is of high dimension, much larger than what others
get. One advantage is that all the variables have a natural
understanding and interpretation in the current approach. This
approach also does not require the Riemann-Roch theorem on $\Jac$ in
order to construct the imbedding. (There is a very ample divisor
lurking in the background, but it is not at all explicit.)

Another advantage of this approach is its elementary nature, since it
only requires rudimentary facts about multi-projective space and the
Segre imbedding. One further drawback of this approach is that it
requires that the base field be extended so that at least three of the
Weierstrass points on $C$ are rational. Over algebraically closed
fields this is not a problem, but some of the other constructions are
not so constrained. 

There are some other things that are not addressed in this
note. First, it is not shown that $\Jac$ is an irreducible variety,
i.e. that it does not have multiple components. It is also not shown
that the equations given for $\Jac$ generate a prime ideal. From a
geometric perspective, it is also not shown here that the defining
equations for $\Jac$ yield a smooth (i.e. non-singular) variety, nor
is it shown that the addition formulas are smooth.

It is however, relatively straightforward (using most computer algebra
systems) to show that $\Jac-\Theta$ is an irreducible affine variety
(i.e. that the defining equations for $\Jac-\Theta$ give a prime
ideal) and that it is smooth (i.e. there are not singular points).

This paper has not been written in the most compact way possible,
and to some it may seem to ramble. However the author feels that it is
important to give an exposition with adequate background for the
construction of genus 2 Jacobians that might be accessible to non-experts,
particularly those with a limited background in modern algebraic geometry.

The general organization of this paper is as follows.  We begin with a
number of remarks about multi-projective varieties, first reviewing
the well known constructions of biprojective varieties and how they can
be turned into projective varieties via the Segre map. This is then
extended in a rather obvious way to multi-projective varieties. Then
there is a short general discussion on divisor classes on
hyperelliptic curves with an explanation of what the $\Theta$-divisor
is and how $\Jac-\Theta$ is readily realizable as an affine
variety. This is applied to the case of genus 2 curves in a very
explicit way to yield a pair of equations in 4 variables that
determine this affine variety.  This is followed by a discussion of
how to add a point of order two, beginning with a special case that
yields a very simple set of formulas for addition. This leads to an
analysis of a biprojective variety that comes fairly close to being
the desired Jacobian, but which still has problems at the projective
closure. In order to patch this problem up, the next idea is to work
in a product of four projective spaces with associated glue equations,
but in order to accomplish this, it is preceeded by a discussion of
how to modify the addition equations for more general Weierstrass
points, and introduces some simplifying notation. This is followed
by an analysis of the infinity types of the associated quadri-projective
variety, which almost gets the right answer, however, there are still
some additional points on this variety that should not be there. We finally
succeed in solving all these problems by going to 8 copies of the affine
variety that are all glued together by 3 different affine Weierstrass points
on the original genus 2 curve.

This results is a rather large set of equations in a high dimensional
multi-projective space, that can be made projective via the Segre
embedding. These equations have a very systematic structure consisting
of a pair of homogeneous equations at each corner of a cube and a set
9 bihomogeneous equations along each edge of the cube, referred to
as glue equations. In the end, it comes out that the points on this
multi-projective variety are in one-to-one correspondence with divisor
classes on the original genus 2 curve.

\head
Multi-Projective Space
\endhead

\subhead
Projective Space
\endsubhead

Working over a fixed base field, $K$, $n$-dimensional projective space,
denoted $\P^n$, is defined by a set of $n+1$ projective coordinates,
$S=\{x_0,x_1,\ldots,x_n\}$, which cannot all simultaneously be 0.  Two
sets of coordinates define the same point in $\P^n$ if all the
coordinates are in a common ratio. Thus $(x_0,x_1,\ldots,x_n)$ and
$(x^\prime_0,x^\prime_1,\ldots,x^\prime_n)$ define the same point in
$\P^n$ if there is some non-zero $\alpha$ such that
$x_i^\prime=\alpha\,x_i$ for $i=0,1,\ldots,n$. 

Affine $n$-space, denoted $\A^n$ corresponds to projective points with
$x_0\ne0$, so $\A^n$ can be understood as just $K^n$, which is a
vector space.  In general, if $V$ is an $n$-dimensional vector space
over $K$, $\P(V)$, which is isomorphic to $\P^n$, can be understood
as the set of lines in $V$.

\subhead
Projective Varieties
\endsubhead

A polynomial, all of whose monomial terms are of degree $d$ in the
variables $S$ is said to be a homogeneous polynomial of degree $d$ in
the variables $S$. If $F(x_0,x_1,\ldots,x_n)$ is homogeneous of degree
$d$, then $$F(\alpha x_0,\alpha x_1,\ldots,\alpha x_n)
=\alpha^nF(x_0,x_1,\ldots,x_n)$$ identically, and therefore the set of
zeros of $F$ give rises to a set of well defined points in $\P^n$. Such
a zero set is a projective variety.  

If $f(x_1,\ldots,x_n)$ is any (not necessarily homogeneous) polynomial
of degree $d$ in the variables $x_1,\ldots,x_n$, there is a corresponding
homogeneous polynomial defined by 
$$\bar f(x_0,x_1,\ldots,x_n)=F(x_0,x_1,\ldots,x_n)=x_0^d\,f(x_1/x_0,\ldots,x_n/x_0).$$
The affine variety defined by $f(x_1,\ldots,x_n)=0$ can be viewed as a
subset of the projective variety defined by $\bar f(x_0,x_1,\ldots,x_n)=0$
corresponding to the affine points, which are defined to be where
$x_0\ne0$. 

Homogeneous polynomials give rise to projective varieties (actually,
projective algebraic sets), and there is a corresponding notion of
projective ideals, which are generated by projective polynomials. If
$I$ is an ideal in $K[x_1,\ldots,x_n]$, then its projective closure
is the projective ideal
$$\bar I = \bigl\{\bar f(x_0,x_1,\ldots,x_n) \vert f(x_1,\ldots,x_n)
\in I\bigr\}$$ 
in $K[x_0,x_1,\ldots,x_n]$ and if $A$ is the affine algebraic set in
$\A^n$ determined by $I$, then its projective closure $\bar A$ is the
projective algebraic set in $\P^n$ determined by $\bar
I$. Computationally, a basis for $\bar I$ can be determined by
homogenizing a Grobner basis for $I$. If $A\subset\A^n$ is an affine
variety, so that $I$ is a prime ideal, then $\bar I$ is also a prime
ideal, and $\bar A$ is a projective variety.

\subhead
Hyperplanes in Projective Space
\endsubhead

A hyperplane in $\A^n$ is defined by a linear equation, i.e.  of the
form 
$$b=\sum_{i=1}^na_ix_i$$ 
with the $a_i$'s not all zero. Any such
hyperplane is isomorphic to a translation of $\A^{n-1}$.

Correspondingly, a hyperplane $\H$ in $\P^n$ is also defined by
a homogeneous linear equation 
$$0=\sum_{i=0}^na_ix_i$$ 
with the $a_i$'s
not all zero. Any such hyperplane is isomorphic to $\P^{n-1}$. The 
intersection of a set of hyperplanes $\cap\H_j$ with each $\H_j$
defined by a homogeneous linear equation $0=\sum_{i=0}^na_{ji}x_i$
is isomorphic to $\P^{k}$ for some $k$.
The dimension $k$ of such a system corresponds to the rank of the matrix
$(a_{ji})$, i.e. the number of independent linear equations.

\subhead
Biprojective Space
\endsubhead

As a set, biprojective $(m,n)$-space is just $\P^m\times\P^n$.  The
Segre map allows $\P^m\times\P^n$ to be defined as a projective
variety in $\P^{(m+1)(n+1)-1}$. If $V_1$ and $V_2$ are vector spaces
over $K$ of dimensions $n$ and $m$ respectively, the Segre map can
be viewed as a natural imbedding of $\P(V_1)\times\P(V_2)$ into
$\P(V_1\otimes V_2)$.

The Segre map
$$\sigma_{n,m}:\, \P^m\times\P^n \rightarrow \P^{(m+1)(n+1)-1}$$
is defined as follows.
If $P_1=(x_0,x_1,\ldots,x_m)\in\P^m$ and
$P_2=(y_0,y_1,\ldots,y_n)\in\P^n$ are given, then by definition,
the Segre map associates the point 
$$(x_0\,y_0, x_1\,y_0, \ldots, x_m\,y_0, x_0\,y_1, x_1\,y_1,
\ldots, x_m\,y_1, x_0\,y_2, \ldots, x_m\,y_m)\in\P^{(m+1)(n+1)-1}$$ 
with $(P_1,P_2)\in\P^m\times\P^n$. 
The image of the Segre map in $\P^{(m+1)(n+1)-1}$ is called the
Segre variety $\Sigma_{m,n}$.   
If the coordinates for $\P^{(m+1)(n+1)-1}$ are $z_k$ for
$k=0,1,\ldots,(m+1)(n+1)-1$, then 
$$z_{j(n+1)+i}=x_i\,y_j$$
defines the Segre map. It is often easier to write
$z_{j(n+1)+i}=z_{i,j}$ for the coordinates of $\P^{(m+1)(n+1)-1}$
where $i$ ranges from $0$ through $m$ and $j$ ranges from $0$ through
$n$.

An important set of relations that hold are
$$z_{i,j}\,z_{k,l}=z_{i,l}\,z_{k,j}$$ 
for $0\le i,k \le m$ and $0\le j,l\le n$, and these relations define
the ideal that determines $\P^m\times\P^n$ as a projective variety in
$\P^{(m+1)(n+1)-1}$. Thus the homogeneous ideal determined by these
quadratic relations determine the Segre variety $\Sigma_{m,n}$.

Note that $i=k$ or $j=l$ give trivial relations
and that by exchanging the order of the terms it may be assumed that
$i<k$ and then by exchanging the left and right sides of the equation,
it may also be assumed that $j<l$, and therefore it is possible to
take $0\le i < k \le m$ and $0\le j < l\le n$ and therefore there
are ${m+1 \choose 2}{n+1 \choose 2}$ such equations that define
$\P^m\times\P^n$ as a projective variety in $\P^{(m+1)(n+1)-1}$. 
However, the difference in dimensions is only $m+n$. This turns
out to be a good example of a projective variety that is not a
complete intersection.

\subhead
Biprojective Varieties
\endsubhead

Suppose that there are now two sets of variables $S_1=\{x_0,x_1,\ldots,x_m\}$ and
$S_2=\{y_0,y_1,\ldots,y_n\}$. A monomial has bidegree $(a,b)$ if the
sum of the degrees of the $x$-variables is $a$ and the sum of the
degrees of the $y$-variables is $b$. Now, a polynomial
$P(x_0,\ldots,x_m,y_0,\ldots,y_n)$ in the $S_1$ and $S_2$ variables is
homogeneous of bidegree $(a,b)$ if all the monomials in $P$ are of
bidegree $(a,b)$. In this case $P(\alpha x_0,\ldots,\alpha x_m,\beta
y_0,\ldots,\beta y_n)=\alpha^a\beta^bP(x_0,\ldots,x_m,y_0,\ldots,y_n)$
identically.  The zero set of such a polynomial gives a well-defined
subset of $\P^m\times\P^n$, and is called a biprojective variety. The
intersection of several such zero sets, each defined by a
bihomogeneous polynomial also gives a biprojective variety.
The Segre map turns any biprojective variety in $\P^m\times\P^n$
into a projective variety in in $\P^{(m+1)(n+1)-1}$.

\subhead
Mapping Bihomogeneous Equations to Homogeneous Equations
\endsubhead

It is often convenient to write $\vec x = (x_0,x_1,\ldots,x_m)$ for
the projective $x$ variables and $\vec y = (y_0,y_1,\ldots,y_m)$ for
the projective $y$ variables Now suppose that $P(\vec x,\vec y)$ is
biprojective of bidegree $(a,b)$, so that $P(\lambda_x\vec
x,\lambda_y\vec y)=\lambda_x^a\,\lambda_y^b\,P(\vec x,\vec y)$ for all
$\lambda_x,\lambda_y\in F$. Suppose that $a<b$ and let $\vec \alpha=
(\alpha_0,\alpha_1,\ldots,\alpha_m)$ be an $m+1$-long vector of non-negative
integers, and define the monomial
$$\vec x^{\vec\alpha}=x_0^{\alpha_0}\,x_1^{\alpha_1}\cdots x_m^{\alpha_m}$$
of degree $\abs{\alpha}=\sum_{i=0}^m\alpha_i$. Note that there are
$m+d \choose d-1$ monomials of degree $d=\abs{alpha}$ in the $m+1$
variables $x_0,x_1,\ldots,x_m$.  With this definition, $\vec
x^{\vec\alpha}\,P(\lambda_x\vec x,\lambda_y\vec y)$ is of bidegree
$(\abs{\alpha}+a,b)$, so in particular, if $\abs{\alpha}=b-a$, then
$\vec x^{\vec\alpha}\,P(\lambda_x\vec x,\lambda_y\vec y)$ is of
bidegree $(b,b)$, i.e. every monomial appearing in $\vec
x^{\vec\alpha}\,P(\lambda_x\vec x,\lambda_y\vec y)$ is a product of
$b$ $x$-variables and $b$ $y$-variables. Now under the Segre map
$$\sigma:\,x_i\,y_j\rightarrow z_{i,j}$$
the polynomial $\vec x^{\vec\alpha}\,P(\lambda_x\vec x,\lambda_y\vec
y)$ maps to a homogeneous polynomial $P^{\alpha}(\vec z)$ of degree
$b$ in the $z$-variables $z_{i,j}$. There are lot of different choices
in how the $x$-variables and the $y$-variables pair up, so strictly
speaking, this is not a completely well defined map if $b>1$.
However, if this is viewed as a subvariety of the Segre variety, all
the relations in the defining equations of the Segre variety make
these different choices irrelevant. In the absence of the Segre
defining equations, it is better to view $P^{\alpha}(\vec z)$ as a
whole set of homogeneous equations in the $z$-variables, all of degree
$b$.

Another thing that can make a difference is the choice of
$\vec\alpha$. Here different choices of $\vec\alpha$ with
$\abs{\vec\alpha}=b-a$ can give different $x$-variables to be combined
with the $y$-variables, and thereby give different $z$-variables.  In
mapping a bihomogenous equation in two set of variables to a
homogenous equation, what needs to be done is to take all possible
choices of $\vec\alpha$ with $\abs{\vec\alpha}=b-a$, which again gives
multiple homogeneous equations associated to each bihomogeneous
equation (unless $b=a$).  Of course, if $b<a$ instead of $a<b$, just
reverse the roles of the two sets of variables, and now multiply
$P(\lambda_x\vec x,\lambda_y\vec y)$ by $\vec x^{\vec\beta}$ where
$\vec\beta=(\beta_0,\beta_1,\ldots,\beta_m)$ is an $n+1$-long vector
of non-negative integers with $\abs{\vec\beta}=\sum_{j=0}^n\beta_j$.

In any case, the set of homogeneous equations that define the
projective variety that is the image of a biprojective variety under
the Segre map, can be determined.

\subhead
Biprojective Hyperplanes
\endsubhead

A biprojective hyperplane $\H$ is now defined by a homogeneous bilinear
equation 
$$0=\sum_{i=0}^m\sum_{j=0}^n a_{i,j}\,x_i\,y_j$$ 
where the $a_{i,j}$ are not all zero. Under the Segre map this becomes
the linear equation $0=\sum_{i=0}^m\sum_{j=0}^n a_{i,j}\,z_{i,j}$.

If $0=\sum_{i=0}^mb_ix_i$ is a hyperplane in only one of the sets ofvariables, 
then this cuts out a subset of $\P^m$ that is isomorphic to
$\P^{m-1}$. There is a corresponding set of hyperplanes in
$\P^{(m+1)(n+1)-1}$ whose intersection give the image of this hyperplane
in the Segre map. The set of bihomogeneous equations is
$0=\sum_{i=0}^mb_ix_iy_j$ for $j=0,1,\ldots,n$ and the corresponding
set of linear hyperplanes in $\P^{(m+1)(n+1)-1}$ is given by the equations
$0=\sum_{i=0}^mb_iz_{i,j}$ for $j=0,1,\ldots,n$. Note that if 
$0=\sum_{i=0}^mb_ix_iy_j$ for $j=0,1,\ldots,n$, then since at least one of
the $y_j$'s must be nonzero, this forces $\sum_{i=0}^mb_ix_i$ to be zero.

\subhead
Infinity Types for Biprojective Varieties
\endsubhead

If $V\subset\P^m\times\P^n$ is a biprojective variety, it is useful to
classify some of the points particularly if
$V\subset\bar{A}\times\bar{B}$ where $A\subset\A^m$ and $B\subset\A^n$
are affine varieties with projective closures $\bar{A}\subset\P^m$ and
$\bar{B}\subset\P^n$, respectively. If $(x_1,\ldots,x_m)$ are affine
coordinates corresponding to projective coordinates
$(x_0,x_1,\ldots,x_m)$ in $\P^m$, then the affine coordinates for a
point in $\A^m$ come from the projective coordinates by simply taking
$x_0=1$, while the projective coordinates that don't come from an
affine point are the result of taking $x_0=0$. Similar ideas apply to
taking the $(y_0,y_1,\ldots,y_n)$ projective coordinates and taking
$y_0=1$ for affine points in $\A^n$ and $y_0=0$ for projective
coordinates in $\P^n$ that don't come from an affine point in
$\A^n$. Thus points in $\P^m\times\P^n$ come in four different
flavors, depending on whether $x_0$ is zero or non-zero and on whether
$y_0$ is zero or non-zero.  This will be referred to as the {\it
infinity type} of a point in $\P^m\times\P^n$.

\subhead
Multi-projective Varieties
\endsubhead

All of this extends very naturally to more than two sets of variables
and to products of more than two projective spaces. If each $S_i$ is a
set of $m_i$ variables and there are $L$ such sets, then a monomial in
the union $S$ of all the $S_i$ can be assigned a multi-degree as an
$L$-tuple of integers $\vec d = (d_1,\ldots,d_L)$ and a polynomial is
homogenous of multi-degree $\vec d$ if each monomial is of multi-degree
$\vec d$. The zero-set of such a homogenous polynomial gives a well
defined subset of the product of projective spaces
$\P^{d_1}\times\cdots\times\P^{d_L}$. The intersection of several such
zero sets of different multi-homogeneous polynomials $f_j$ in $S$, the
set of all variables, is the set of points annihilated by all elements
of the multi-graded ideal generated by the $f_j$.

There is a natural generalization of the Segre map denoted
$$\sigma_{\vec d}:\,\P^{d_1}\times\cdots\times\P^{d_L}\rightarrow\P^D$$
where $D=\prod_{i=1}^L(d_i+1) - 1$, and the image of
$\P^{d_1}\times\cdots\times\P^{d_L}$ in $\P^D$ is the generalized Segre
variety, denoted $\Sigma_{\vec d}$.

If $S_i=\{x_{i,0},x_{i,1},\ldots,x_{i,d_i}\}$ is the $i$-th set of
variables, representing a point in $\P^{d_i}$, and if $\vec x_i = 
(x_{i,0},x_{i,1},\ldots,x_{i,d_i})$ is the corresponding vector
of variables, then
$$\sigma_{\vec d}: x_{1,k_1}\,x_{2,k_2}\cdots x_{d,k_d}
\rightarrow z_{k_1,k_2,\ldots,k_d}$$
is the mapping from a product of projective spaces to a new projective
space in the $z_{\vec k}$-coordinates, with corresponding maps in the
$z_{\vec k}$-variables. The way to go from multi-homogeneous defining
equations to homogeneous defining equations is by multiplying by
appropriate ``slack'' monomials, just as in the case of bihomogeneous
polynomials. Analogous concepts apply to multi-projective hyperplanes.

Infinity types for points in multi-projective spaces are defined in analogy
to their concept in biprojective space. Starting with affine varieties
$A_i\subset\A^{d_i}$ with projective closures $\bar A_i\subset\P^{d_i}$,
points in $\bar A_1 \times \bar A_2 \times\cdots\times \bar A_L
\subset \P^{d_1}\times\P^{d_2}\times\cdots\P^{d_L}$ are of infinity
type $\bigl(x_{1,0}=c_1,x_{2,0}=c_2,\ldots,x_{L,0}=c_L\bigr)$ with
$c_1,c_2,\ldots,c_L\in\{0,1\}$ depending on whether the coordinates
are in the affine or non-affine part of the corresponding component.
This is a very useful idea when glueing different parts of projective
varieties together along common affine subsets, and then looking at
what happens on the different parts of the projective closures.

\head
Hyperelliptic Curves and the Picard Group
\endhead

Let $F$ be a field that is not of characteristic 2.
Let $f(x)$ be a monic polynomial over the base field $F$ 
of degree $2\,g+1$ with no repeated roots.
The curve $\Caff$ defined by the affine equation
$$\Caff:\,y^2=f(x)$$ 
is hyperelliptic of genus $g$ with hyperelliptic involution
$\iota:\,\Caff\rightarrow\Caff$ defined by $\iota:\,(x,y)
\mapsto(x,-y)$.  Its projective closure $C$ is defined by the
homogenous equation
$$C:\, Y^2\,Z^{2g-1}=Z^{2g+1}\,f(X/Z)$$
and is obtained by adding the point $P_\infty$
at infinity with projective coordinates $(1,0,0)$
to $\Caff$ where the standard inclusion $\Caff\rightarrow C$
is given by $(x,y)\mapsto(x,y,1)$ in projective coordinates.
The hyperelliptic involution is extended by defining
$\iota(P_\infty)=P_\infty$.

A Weierstrass point on $C$ is a point $P\in C$ such that
$\iota(P)=P$.  In addition to $P_\infty$, the affine Weierstrass
points of $C$ are of the form $(\rho,0)$ where $\rho$ is a
root of $f(x)$. Thus, there are exactly $2\,g+2$ Weierstrass
points on $C$ if $C$ is of genus $g$.

\subhead
Divisors and Divisor Classes on Hyperelliptic Curves
\endsubhead

A divisor $D$ on $C$ is just a finite formal sum of points with
multiplicities, i.e. $D=\sum_{i=1}^k n_i\,P_i$ where $P_i\in C$.
The degree of $D$ is defined by $\deg(D)=\sum_{i=1}^k n_i$. Every
non-zero function $h$ on $C$ has an associated divisor
$$\div(h)=\sum_{P\in C} \deg_P(h)\cdot P$$
and since every function on $C$ has only finitely many zeros and poles,
this is actually a finite sum. The degree of the divisor of any non-zero
function on $C$ is 0.

Two divisors $D$ and $D^\prime$ are equivalent if their difference
is the divisor of a function, i.e. $D\equiv D^\prime$ if there
is some function $h$ on $C$ such that $D=D^\prime+\div(h)$. The
group of equivalence classes of divisors on $C$ is the Picard group
$\hbox{Pic}(C)$ and the group of equivalence classes of degree 0
is denoted $\Pic$.

The geometric realization of $\Pic$ is the Jacobian variety $\Jac$,
which is a projective algebraic variety of dimension $g$.

A divisor $D$ of degree 0 is semi-reduced if $D$ is of the form
$$D=\sum_{i=1}^k P_i - k\cdot P_\infty$$
where the $P_i$ are all on $\Caff$ and are such that if 
$P_i=\iota(P_j)$ then $i=j$, i.e. no affine point and its
hyperelliptic involute and appear in the support of $D$,
and if an affine Weierstrass point is in the support of $D$ then
it appears with multiplicity 1. A semi-reduced divisor is reduced
if $k\le g$.

\proclaim{Theorem} Every divisor class of degree 0 on $C$ contains
exactly one reduced divisor.
\endproclaim

This follows from the Riemann-Roch theorem, but can also be
proved by extending the theory of reduction of quadratic
forms to polynomials over the base field $F$. The composition
of quadratic forms corresponds to addition of divisor classes.

\subhead
Polynomial Representation of Reduced and Semi-Reduced Divisors
\endsubhead

If $D=\sum_{i=1}^k P_i - k\cdot P_\infty$ is a semi-reduced divisor
there is a polynomial $U(x;D)$ of degree $k$ whose roots are the
$x$-coordinates of the $P_i$, i.e.
$$U(x;D)=\prod_{i=1}^k \bigl(x-x(P_i)\bigr)$$
and there is an associated polynomial $V(x;D)$ of degree at most
$k-1$ that interpolates the $y$-coordinates, i.e.
$$V(x(P_i);D)=y(P_i)$$ 
for $i=1,\ldots,k$. Note that if all the points in the support of $D$
are distinct, then $f(x)-V(x;D)^2$ is 0 for $x=x(P_i)$, and therefore
$U(x;D)$ divides $f(x)-V(x;D)^2$. If $D$ has points of multiplicity greater
than 1, then this is the defining property of $V(x;D)$ along with
the requirement $\deg V(x;D) < k$. Thus there is a unique polynomial
$W(x;D)$ such that 
$$f(x)-V(x;D)^2=U(x;D)\,W(x;D).$$
Furthermore, and any such triple $\bigl(U(x),V(x),W(x)\bigr)$
with $f=V^2+U\,W$ defines a unique semi-reduced divisor. 
The representation theorem above simply states that every divisor
class of degree 0 has a unique $\bigl(U(x),V(x),W(x)\bigr)$ with
$\deg U(x) \le g$ and $\deg V(x) < g$.

The set of divisor classes represented by a reduced divisor
$D=\sum_{i=1}^k P_i - k\cdot P_\infty$ with $k < g$ is denoted
$\Theta$. Geometrically $\Theta$ is a subvariety of $\Jac$ of
codimension 1, and is referred to as the $\Theta$-divisor.
(This can be confusing, since $\Theta$ is not a divisor on $C$,
but rather on $\Jac$.)

\subhead
$\Jac-\Theta$ as an affine variety
\endsubhead

$\Jac-\Theta$ can be given the structure of an affine variety
(actually, just an affine algebraic set) as follows:
Write 
$$\eqalign{
U(x)&=x^{g}+\sum_{i=0}^{g-1} u_i\,x^i \cr
V(x)&=\sum_{i=0}^{g-1} u_i\,x^i \cr
W(x)&=x^{g+1}+\sum_{i=0}^{g} u_i\,x^i \cr
}$$
for variables 
$\bolds=\bigl\{u_0,\ldots,u_{g-1},v_0,\ldots,v_{g-1},w_0,\ldots,w_{g}\bigr\}$
and write 
$$f(x)=x^{2g+1}+\sum_{i=0}^{2g}a_i\,x^i$$
and now simply equate coefficients in $f=V(x)^2+U(x)\,W(x)$. These
equations on the coefficients define $\Jac-\Theta$ as an affine variety.
Alternatively, the $w_i$'s can be eliminated by simply reducing $f(x)-V(x)^2$
modulo $U(x)$ generically, obtaining a polynomial of degree $g-1$ in $x$
with coefficients that are polynomials in the $u_i$'s and $v_i$'s, and
the variety now follows by requiring that these coefficients all be 0.

\head
Genus 2 curves
\endhead

Let $f(x)$ be a monic quintic polynomial over the base field $F$ having 
no multiple roots and write
$$\eqalign{f(x)
&=x^5+a_4\,x^4+a_3\,x^3+a_2\,x^2+a_1\,x+a_0 \cr
&=\bigl(x-\rho^{(1)}\bigr)\,\bigl(x-\rho^{(2)}\bigr)\,\bigl(x-\rho^{(3)}\bigr)
 \,\bigl(x-\rho^{(4)}\bigr)\,\bigl(x-\rho^{(5)}\bigr) \cr
}$$
with $\rho^{(i)}\ne\rho^{(j)}$ for $i \ne j$. These roots $\rho^{i}$ may be in
finite extension of the base field $F$. It is also be useful to expand
$f(x)$ around the $\rho{(i)}$'s as
$$f(x)
=\bigl(x-\rho^{(1)}\bigr)^5+a_4^{(i)}\,\bigl(x-\rho^{(1)}\bigr)^4
+a_3^{(i)}\,\bigl(x-\rho^{(1)}\bigr)^3+a_2^{(i)}\,\bigl(x-\rho^{(1)}\bigr)^2
+a_1^{(i)}\,\bigl(x-\rho^{(1)}\bigr)
$$ 
with $a_1^{(i)}\ne0$ for $i=1,\ldots,5$. For compactness of
notation. it is useful to write 
$$\eqalign{\
\bolda&=(a_0,a_1,a_2,a_3,a_4,1) \cr
\boldx&=(1,x,x^2,x^3,x^4,x^5) \cr
}$$ 
so $f(x)=\bolda\cdot\boldx$. It is also useful to write
$$\eqalign{
\bolda^{(i)}& =\bigl(0,a_1^{(i)},a_2^{(i)},a_3^{(i)},a_4^{(i)},1\bigr) \cr
\boldx^{(i)}&=\bigl(1,(x-\rho^{(i)}),(x-\rho^{(i)})^2,(x-\rho^{(i)})^3,(x-\rho^{(i)})^4,(x-\rho^{(i)})^5\bigr) \cr}$$
so $f(x)=\bolda^{(i)}\cdot\boldx^{(i)}$, as well for $i=1,\ldots,5$.

\subhead
Divisor classes on genus 2 curves
\endsubhead

In the particular case of genus 2, a general degree 0 divisor class
will be represented by a reduced divisor of the form
$$D=P_1+P_2-2\,P_\infty$$ 
where $P_1$ and $P_2$ are affine points on $C$ and $\iota(P_1)\ne P_2$
(although $P_1=P_2$ is allowed if $P_1$ is not a Weierstrass
point). In this case the $\Theta$-divisor consists of those divisor classes
on $C$ represented by reduced divisors of the form
$$D=P-P_\infty$$
for $P$ an affine point on $C$, along with the divisor $0$. Alternatively,
$\Theta$ consists of those divisor classes with reduced divisor $P-P_\infty$
as $P$ ranges over all of $C$. Thus $\Theta$ is essentially an image of $C$
imbedded in $\Jac$ under the Abel-Jacobi mapping.

There are basically three types of degree 0 divisor classes; they are
represented by reduced divisors of one of the three following types:
\vskip1pt $(i)$ $P_1+P_2-2\,P_\infty$ where $P_1,P_2\in\Caff$ with $P_1\ne\iota(P_2)$;
\vskip1pt $(ii)$ $P_1-P_\infty$ where $P_1\in\Caff$;
\vskip1pt $(iii)$ $0$.
\vskip2pt\noindent 
The divisors of type $(ii)$ and $(iii)$ constitute the set $\Theta$
known as the theta-divisor, and represent a copy of $C$ inside of
$\Jac$.  Notionally, $\Theta$ is a one-dimensional object, while
$\Jac$ is a two-dimensional object (and more generally, on
hyperelliptic curves of genus $g$, $\Theta$ is a codimension one
subvariety of the $g$-dimensional Jacobian variety).

\subhead
Points of order 2 on $\Jac$
\endsubhead

In general, an abelian variety of dimension $g$ has $2^{2g}$ points of
order two, so in the present case there are 16 points of order 2 on $\Jac$.
If $P_i=(\rho_i,0)$ is an affine Weierstrass point on $C$, then
$X_i=P_i-P_\infty$ is a reduced divisor representing a point of order 2 on
$\Jac$ since the divisor of the function $x-\rho_i$ is
$2\,P_i-2\,P_\infty$. There are 5 such points, all on $\Theta$, and also
there is the divisor $0$ on $\Theta$, making a total of 6 points on $\Theta$
of order 2. Furthermore, if $P_i$ and $P_j$ are distinct
affine Weierstrass points on $C$, then $X_i+X_j=P_i+P_j-2\,P_\infty$ is also a
reduced divisor representing a point of order 2 on $\Jac$. This accounts
for all 16 points on $\Jac$ of order 2. It is also worth noting that the
divisor of the function $y$ on $C$ is $P_1+P_2+P_3+P_4+P_5-5\,P_\infty$,
i.e $X_1+X_2+X_3+X_4+X_5=0$ on $\Jac$.

\subhead
Defining equations for $\Jac-\Theta$ as an affine variety in genus 2
\endsubhead

Let $\bolds=\{u_0,u_1,v_0,v_1\}$ be a set of (affine) variables
and let
$$\eqalign{
U(x) &= x^2+u_1\,x+u_0 \cr
V(x) &= v_1\,x+v_0 \cr
}$$
be a pair of polynomials that represents a reduces divisor
on $\Jac-\Theta$, so in particular they satisfy the relationship
$$U(x) \vert f(x)-V(x)^2$$
and it is instructive to write
$$f(x) - V(x)^2 = e_1\,x+e_0 \bmod U(x)$$
with $e_1=e_1(u_0,u_1,v_0,v_1)$ and $e_0=e_0(u_0,u_1,v_0,v_1)$
for a pair of polynomials $e_0,e_1\in F[u_0,u_1,v_0,v_1]$. 
Specifically, these polynomials are
$$\eqalign{
&e_0(v_1,v_0,u_1,u_0) \cr 
&=\! v_1^2\,u_0 \!-\! v_0^2 \!+\! u_1^3\,u_0 \!-\! a_4\,u_1^2\,u_0 \!-\! 2\,u_1\,u_0^2 \!+\! a_3\,u_1\,u_0 \!+\!
    a_4\,u_0^2 \!-\! a_2\,u_0 \!+\! a_0 \cr
&e_1(v_1,v_0,u_1,u_0) \cr 
&=\! v_1^2\,u_1 \!-\! 2\,v_1\,v_0 \!+\! u_1^4 \!-\! a_4\,u_1^3 \!-\! 3\,u_1^2\,u_0 \!+\! a_3\,u_1^2 \!+\!
    2\,a_4\,u_1\,u_0 \!-\! a_2\,u_1 \!+\! u_0^2 \!-\! a_3\,u_0 \!+\! a_1 \cr
}$$
Any point $(u_0,u_1,v_0,v_1)\in F^4$ where these two polynomials
vanish simultaneously gives a pair of polynomials
$\bigl(U(x),V(x)\bigr)
=\bigl(U(x;u_0,u_1,v_0,v_1),V(x;u_0,u_1,v_0,v_1)\bigr)$ with $U(x)$
monic, such that $U(x) \vert f(x)-V(x)^2$. 

\proclaim{Proposition}
Let $\bolds=(u_0,u_1,v_0,v_1)$ be a sequence of coordinates in $F^4$,
and let $A(\bolds)\in K^4$ be the zero set of the following
pair of polynomials:
$$\eqalign{
e_0&: a_0 - a_2 u_0 + a_4 u_0^2 + a_3 u_0 u_1 - 2 u_0^2 u_1 - a_4 u_0 u_1^2 + u_0 u_1^3 - v_0^2 + u_0 v_1^2 \cr
e_1&: a_1 - a_3 u_0 + u_0^2 - a_2 u_1 + 2 a_4 u_0 u_1 + a_3 u_1^2 - 3 u_0 u_1^2 - a_4 u_1^3 + u_1^4 - 2 v_0 v_1 + u_1 v_1^2 \cr
}$$
Then the points on the affine algebraic set $A(\bolds)$ are in one-to-one
correspondence with elements of $\Jac-\Theta$.
\endproclaim

These two polynomials define a
prime ideal in $F[u_0,u_1,v_0,v_1]$ whose associated affine variety is
$\Jac-\Theta$.

\head
Addition of a Point of Order 2: The Case of a Distinguished Weierstrass Point
\footnote{In this section only, the equations are modified from the more general
equations above.}
\endhead

It is useful to consider the case of $a_0=0$ in the equation for $C$,
which corresponds to $C$ having a Weierstrass point at $P_0=(0,0)$. 
The equations for $A(\bolds)$
are modified by setting $a_0=0$ in the polynomials $e_0$ and $e_1$. However,
$a_0$ appears only in equation $e_0$, which now becomes
$$
e_0 =  - a_2 u_0 + a_4 u_0^2 + a_3 u_0 u_1 - 2 u_0^2 u_1 - a_4 u_0 u_1^2 + u_0 u_1^3 - v_0^2 + u_0 v_1^2 .$$

Setting $u_0=0$ in the equations for $e_0$ and $e_1$ gives
the pair of equations
$$\eqalign{
e_0: 0 &= v_0^2 + u_0 v_1^2 \cr
e_1: 0 &= a_1 - a_2 u_1 + a_3 u_1^2 - a_4 u_1^3 + u_1^4 - 2 v_0 v_1 + u_1 v_1^2 \cr
}$$
so $v_0=0$ and the second equation becomes
$$0=a_1 - a_2 u_1 + a_3 u_1^2 - a_4 u_1^3 + u_1^4 - 2 v_0 v_1 + u_1 v_1^2$$
just as before.

\proclaim{Proposition} 
Suppose that $0$ is a root of $f(x)$ so $P_0=(0,0)$ is an affine
Weierstrass point on $C$. Let ${\boldkey P}=(u_0,u_1,v_0,v_1)$ be the
coordinates of a point on $A(\bolds)$.  If $u_0=0$ then $v_0=0$ and
$P=(-u_1,u_1v_1)$ is an affine point on $C$ with $P\ne P_0$, and
${\boldkey P}$ corresponds to the reduced divisor $P+P_0-2\,P_\infty$
as an element of $\Jac-\Theta$.  \endproclaim

\subhead
Adding a point of order 2
\endsubhead

The general procedure for adding two elements of $\Pic$ 
(but which are assumed to all have distinct points in their support) that are
represented by polynomial pairs $\bigl(U^{(1)}(x),V^{(1)}(x)\bigr)$ and
$\bigl(U^{(2)}(x),V^{(2)}(x)\bigr)$ (which need not lie outside of $\Theta$)
is as follows: 
\vskip 2pt
(1) Find a polynomial $\breve V(x)$ that is of degree at most 3 and
such that $$\breve V(x) \equiv V^{(i)}(x) \bmod U^{(i)}(x)$$ for
$i=1,2$.
\vskip 2pt
(2) Find the points where $y=\breve V(x)$ intersects the curve
$y^2=f(x)$.  This amounts to solving $f(x) - \breve V(x)^2=0$.  Some
of these will be points already represented by the original polynomial
pairs.  The points that are not are the ones whose involutes are in
the sum of the divisors.  Thus $$U^{(3)}(x) = \left(f(x) -
\breve V(x)^2\right)/U^{(1)}(x)U^{(2)}(x)$$
where $U^{(3)}(x)$ has as roots the $x$-coordinates of the desired
points.
\vskip 2pt
(3) To find the interpolating polynomial for the $y$-coordinates, just
note that $$\breve V(x) \equiv -V^{(3)}(x) \bmod U^{(3)}(x)$$
since the points on the intersection must be involuted to get the
points in the sum.
\vskip 2pt\noindent
This procedure works perfectly well as long as there are no common
factors between the $U^{(i)}(x)$'s.

This procedure is particularly simple to work out in one very
interesting case, namely, when one of the points on $\Pic$ being added
is the distinguished 2-division point 
$$X_0 = \hbox{Cl}(P_0-P_\infty).$$
For this point the representing polynomial pair is
$\left(U(x;X_0)\,V(x;X_0)\right) = (x,0)$.  
For simplicity of notation, let $X\in\Pic$ be arbitrary except for the 
requirement that  $X,X+X_0\not\in\Theta$, and set
$$\eqalign{
U(x) &= U(x;X) \cr
V(x) &= V(x;X) \cr
\hat U(x) &= U(x;X + X_0) \cr
\hat V(x) &= V(x;X + X_0) \cr
}$$
and in terms of explicit coefficients, write
$$\eqalign{
U(x) &= x^2 + u_1\,x + u_0 \cr
V(x) &= v_1\,x + v_0 \cr
\hat U(x) &= x^2 + \uhat_1\,x + \uhat_0 \cr
\hat V(x) &= \vhat_1\,x + \vhat_0 \cr
}$$
with a goal of finding relations between the two sets of coordinates
$\bolds\!=\!\left(u_0, u_1, v_0, v_1\right)$ and $\boldhats=\left(\uhat_0,
\uhat_1, \hat v_0, \vhat_1\right)$.

In the case at hand, the congruence relation on $\breve V(x)$ imposed by
$X_0$ is that $x\vert \breve V$.  This means that it is possible to
write 
$$\breve V(x) = \tilde V(x)\,x$$
for some polynomial $\tilde V(x)$.  Then 
$$f(x) - \breve V(x)^2 = U(x)\,\hat U(x)\,x$$
so dividing by $x$ gives
$$x^4+a_4\,x^3+a_3\,x^2+a_2\,x+a_1 =  \tilde V(x)^2\,x + U(x)\,\hat U(x)$$
and then setting $x=0$ gives the relation
$$a_1 = u_0\,\uhat_0.$$
Since $a_1\ne 0$, this shows that neither $u_0$ nor $\uhat_0$ can be
0 if neither $X$ nor $X+X_0$ is on $\Theta$.  This corresponds to what
can be seen by looking at reduced divisors directly, and gives an
indication of how the algebraic equations reflect the group law.

The other congruence relations on
$\breve V$ are then satisfied by writing
$$\eqalign{
\breve V(x) &= V(x) - {v_0 \over u_0} U(x) \cr
\breve V(x) &= - \hat V(x) + {\vhat_0 \over \uhat_0} \hat U(x) \cr
}$$
which now gives
$$V(x)+\hat V(x) = {\vhat_0 \over \uhat_0} \hat U(x) + {v_0 \over u_0} U(x)
=a_1^{-1}\left(\vhat_0\,u_0\,\hat U(x) + v_0\,\uhat_0\,U(x)\right)$$
in view of the relation between $u_0$ and $\uhat_0$.  Equating the
coefficients of $x^2$ and $x$ and using the fact that $a_1 \ne 0$ now gives
$$\eqalign{
0 &= \vhat_0\,u_0 + v_0\,\uhat_0 \cr
a_1\left(v_1+\vhat_1\right) &= \vhat_0\,u_0\,\uhat_1 + v_0\,\uhat_0\,u_1 \cr
}$$
as further relations between the variables.
Substituting the first equation into the second gives
$$a_1\left(v_1+\vhat_1\right) = \vhat_0\,u_0\,\uhat_1 - \vhat_0\,u_0\,u_1$$
and multiplying by $\uhat_0$ and dividing by $a_1$ then gives
$$\uhat_0\left(v_1+\vhat_1\right) 
= \vhat_0\left(\uhat_1 - u_1\right)$$
and there is also the equation
$$u_0\left(v_1+\vhat_1\right) 
= v_0\left(u_1 - \uhat_1\right)$$
where the roles of the hatted and the unhatted variables are reversed.

The above pair of equations for $\breve V(x)$ also shows that
$\breve V(x)$ is of degree at most 2, which means that $\tilde V(x)$
is of degree at most 1.  The coefficient of $x^2$ in $\breve V(x)$ is
seen to be $-{v_0 / \uhat_0}$ and also ${\vhat_0 / u_0+}$
and therefore the coefficient of $x^3$ in $\tilde V(x)^2$ is $-{v_0\hat
v_0 / a_1}$.  Therefore equating coefficients of $x^3$ in the
equation above gives
$$a_1\,a_4 = -v_0\,\vhat_0 + a_1\,u_1 + a_1\,\uhat_1$$
as another equation.

To summarize so far, the following set of equations holds between the two sets
of coordinates:
$$\eqalign{
0 = g_1(\bolds,\boldhats) &= u_0\,\uhat_0 - a_1 \cr
0 = g_2(\bolds,\boldhats) &= \vhat_0\,u_0 + v_0\,\uhat_0 \cr
0 = g_3(\bolds,\boldhats) &= v_0\,\vhat_0 + a_1\left(a_4 - u_1 - \uhat_1\right) \cr
0 = g_4(\bolds,\boldhats) &= \uhat_0\left(v_1+\vhat_1\right) + \vhat_0\left(u_1 - \uhat_1\right) \cr
0 = g_5(\bolds,\boldhats) &= u_0\left(v_1+\vhat_1\right) + v_0\left(\uhat_1 - u_1\right) \cr
}$$
where for convenience, the polynomials have been labelled
$g_1,g_2,g_3,g_4,g_5$ for ease of reference.
These equations suffice to determine the values of $\boldhats=\left(\uhat_0,
\uhat_1, \vhat_0, \vhat_1\right)$ in terms of the values of
$\bolds=\left(u_0, u_1, v_0, v_1\right)$ as long as $u_0 \ne 0$, and vice
versa.

\subhead
Additional relations for adding a distinguished point of order 2
\endsubhead

These are not the only relations that are of interest, however.  There
are other interesting coordinate equations that can be seen by looking
at the projective closure of the graph of the two affine varieties
that are related by the addition of a point of order 2.  
One of the concerns that naturally arise from looking at the equations
for $e_0$ and $e_1$ is that they are fairly high degree. 

The goal will now be to derive some relations between hatted and
unhatted variables that are of lower degree. The possibility of doing
this is suggested by some of the degree lowering relations that are
apparent in some of the glue equations (e.g. $u_0\,\uhat_0=a_1$),
applying this relation to the $e_0$ and $e_1$ polynomials.  For now,
just note that these new relations depend on $u_0$ and $\uhat_0$ both
being non-zero.

Start by noting that
$$\eqalign{
e_0&(\bolds,\boldhats) \uhat_0 + g_2(\bolds,\boldhats) v_0 \cr
&= u_0 \left(a_4\,\uhat_0\,u_0-a_2\,\uhat_0  + a_3\,\uhat_0\,u_1 -
2\,\uhat_0\,u_0\,u_1 - a_4\,\uhat_0\,u_1^2 + \uhat_0\,u_1^3 + \hat
v_0\,v_0 + \uhat_0\,v_1^2\right)\cr}$$
and therefore, assuming $u_0 \ne 0$ this gives the relation
$$0 = h_0(\bolds,\boldhats) = 
-a_2\,\uhat_0 + a_4\,\uhat_0\,u_0 + a_3\,\uhat_0\,u_1 -
2\,\uhat_0\,u_0\,u_1 - a_4\,\uhat_0\,u_1^2 + \uhat_0\,u_1^3 + \hat
v_0\,v_0 + \uhat_0\,v_1^2$$
where the polynomial $h_0$ has been defined here for convenience.
Next define
$$\eqalign{h_1(\bolds,\boldhats) 
&= e_1(\bolds,\boldhats)\,\uhat_0 - h_0(\bolds,\boldhats)\,u_1 \cr
&= a_1\,\uhat_0 - a_3\,\uhat_0\,u_0 + \uhat_0\,u_0^2 + a_4\,\hat
u_0\,u_0\,u_1 - \uhat_0\,u_0\,u_1^2 - u_1\,\vhat_0\,v_0 - 2\,\hat
u_0\,v_0\,v_1 \cr
}$$
and now use the relation $u_0\,\uhat_0 = a_1$ to get that $h_2(\bolds,\boldhats)=0$
where
$$h_2(\bolds,\boldhats) = 
  -a_1\,a_3 + a_1\,\uhat_0 + a_1\,u_0 + a_1\,a_4\,u_1 - a_1\,u_1^2 -
u_1\,\vhat_0\,v_0 - 2\,\uhat_0\,v_0\,v_1$$
and finally
$$\eqalign{g_6(\bolds,\boldhats) 
&= h_2(\bolds,\boldhats) - g_3(\bolds,\boldhats)\,u_1 \cr
&=  -a_1\,a_3 + a_1\,\uhat_0 + a_1\,u_0 + a_1\,\uhat_1\,u_1 -
2\,u_1\,\vhat_0\,v_0 - 2\,\uhat_0\,v_0\,v_1 \cr
}$$
must also be a relation between the unhatted and the hatted variables.
There must also be a relation with the roles of the hatted and
unhatted variables reversed, i.e.
$$g_7(\bolds,\boldhats) =   
  -a_1\,a_3 + a_1\,\uhat_0 + a_1\,u_0 + a_1\,\uhat_1\,u_1 -
2\,\uhat_1\,\vhat_0\,v_0 - 2\,u_0\,\vhat_0\,\vhat_1
$$
must also be 0.

Next, note that 
$$0=g_3(\bolds,\boldhats)\,u_0^2-v_0\,u_0\,g_2(\bolds,\boldhats)+v_0^2\,g_1(\bolds,\boldhats)
=a_1\,\left(a_4\,u_0^2 - \uhat_1\,u_0^2 - u_0^2\,u_1 - v_0^2\right)$$
and therefore 
$$0 = i_0(\bolds,\boldhats) 
= a_4\,u_0^2 - \uhat_1\,u_0^2 - u_0^2\,u_1 - v_0^2$$
since $a_1 \ne 0$.  Then,
$$0= e_0(\bolds,\boldhats)-i_0(\bolds,\boldhats) 
= u_0\,\left(-a_2 + \uhat_1\,u_0 + a_3\,u_1 - u_0\,u_1 - a_4\,u_1^2 +
u_1^3 + v_1^2\right)$$
and so
$$0 = i_1(\bolds,\boldhats) 
= -a_2 + \uhat_1\,u_0 + a_3\,u_1 - u_0\,u_1 - a_4\,u_1^2 + u_1^3 + v_1^2$$
since $u_0 \ne 0$.  Next, set
$$\eqalign{
0 &= i_2(\bolds,\boldhats) \cr 
&= \vhat_0\,i_1(\bolds,\boldhats) + \left(u_1 - \uhat_1\right)\,g_2(\bolds,\boldhats) + g_4\,\uhat_0
- \left(v_1 + \vhat_1\right)\,g_1(\bolds,\boldhats) \cr &= 
  -a_2\,\vhat_0 + a_3\,u_1\,\vhat_0 - a_4\,u_1^2\,\vhat_0 + u_1^3\,\hat
v_0 + a_1\,\vhat_1 + a_1\,v_1 + \vhat_0\,v_1^2\cr
}$$ 
and finally set
$$\eqalign{
0&= g_8(\bolds,\boldhats) 
= i_2(\bolds,\boldhats) u_1 - \vhat_0 e_1(\bolds,\boldhats) \cr 
&=  -a_1 \vhat_0 + a_3 u_0 \vhat_0 - u_0^2 \vhat_0 -
2 a_4 u_0 u_1 \vhat_0 + 3 u_0 u_1^2 \vhat_0 + a_1 u_1 \vhat_1 +
a_1 u_1 v_1 + 2 \vhat_0 v_0 v_1 \cr
}$$
as another equation relating the hatted and unhatted variables.  In
addition, there is equation
$$\eqalign{
0&= g_9(\bolds,\boldhats) \cr
&=  -a_1 v_0 + a_3 \uhat_0 v_0 - \uhat_0^2 v_0 -
2 a_4 \uhat_0 \uhat_1 v_0 + 3 \uhat_0 \uhat_1^2 v_0 +
a_1 \uhat_1 v_1 + a_1 \uhat_1 \vhat_1 + 2 v_0 \vhat_0 \vhat_1 
}$$
obtained by reversing the roles of the hatted and the unhatted variables.

To summarize, the following equations hold between the coordinates for
$X$ and $X+X_0$:
$$\eqalign{
0 = g_1(\bolds,\boldhats) &= u_0\,\uhat_0 - a_1 \cr
0 = g_2(\bolds,\boldhats) &= \vhat_0\,u_0 + v_0\,\uhat_0 \cr
0 = g_3(\bolds,\boldhats) &= v_0\,\vhat_0 + a_1\left(a_4 - u_1 - \uhat_1\right) \cr
0 = g_4(\bolds,\boldhats) &= \uhat_0\left(v_1+\vhat_1\right) + \vhat_0\left(u_1 - \uhat_1\right) \cr
0 = g_5(\bolds,\boldhats) &= u_0\left(v_1+\vhat_1\right) + v_0\left(\uhat_1 - u_1\right) \cr
0 = g_6(\bolds,\boldhats) 
   &= -a_1\,a_3 + a_1\,\uhat_0 + a_1\,u_0 + a_1\,\uhat_1\,u_1  -
2\,u_1\,\vhat_0\,v_0 - 2\,\uhat_0\,v_0\,v_1 \cr
0 = g_7(\bolds,\boldhats) 
   &= -a_1\,a_3 + a_1\,\uhat_0 + a_1\,u_0 + a_1\,\uhat_1\,u_1  -
2\,\uhat_1\,\vhat_0\,v_0 - 2\,u_0\,\vhat_0\,\vhat_1 \cr
0 = g_8(\bolds,\boldhats)
   &= -a_1\,\vhat_0 + a_3\,u_0\,\vhat_0 - u_0^2\,\vhat_0 -
2\,a_4\,u_0\,u_1\,\vhat_0 + 3\,u_0\,u_1^2\,\vhat_0 + a_1\,u_1\,\vhat_1 \nl +
a_1\,u_1\,v_1 + 2\,\vhat_0\,v_0\,v_1 \cr
0 = g_9(\bolds,\boldhats)
   &= -a_1\,v_0 + a_3\,\uhat_0\,v_0 - \uhat_0^2\,v_0 -
2\,a_4\,\uhat_0\,\uhat_1\,v_0 + 3\,\uhat_0\,\uhat_1^2\,v_0 +
a_1\,\uhat_1\,v_1 \nl + a_1\,\uhat_1\,\vhat_1 + 2\,v_0\,\vhat_0\,\vhat_1 \cr
}$$
and it is worth noting that $g_1$, $g_2$, and $g_3$ are symmetric,
i.e. that $g_1(\bolds,\boldhats)=g_1(\boldhats,\bolds)$,
$g_2(\bolds,\boldhats)=g_2(\boldhats,\bolds)$, and
$g_3(\bolds,\boldhats)=g_3(\boldhats,\bolds)$, while $(g_4,g_5)$,
$(g_6,g_7)$, and $(g_8,g_9)$ are complementary pairs, i.e. 
$g_5(\bolds,\boldhats)=g_4(\boldhats,\bolds)$, 
$g_7(\bolds,\boldhats)=g_6(\boldhats,\bolds)$, and
$g_9(\bolds,\boldhats)=g_8(\boldhats,\bolds)$. It is also worth
noting that the bidegrees of $g_1$, $g_2$ and $g_3$ are $(1,1)$,
the bidegrees of $g_4$ and $g_7$ are $(1,2)$, the bidegrees of
$g_5$ and $g_6$ are $(2,1)$, and the bidegrees of $g_8$ and $g_9$
are $(3,1)$ and $(1,3)$, respectively.

\proclaim{Proposition} 
Supppose $f(x)=x^5+a_4\,x^4+a_3\,x^3+a_2\,x^2+a_1\,x$ has no repeated
roots and let $C$ be the curve $y^2=f(x)$. Let $X_0$ be the point on
$\Jac$ corresponding to the divisor class $(0,0)-P_\infty$. Suppose
that $X\in\Jac-\Theta$ and that $X+X_0\in\Jac-\Theta$ also. Let
$X$ be represented by a quadruple $\bolds=(u_0,u_1,v_0,v_1)$
corresponding to polynomials $U(x)=x^2+u_1\,x+u_0$ and
$V(x)=v_1\,x+v_0$ and such that $U(x) \vert f(x)-V(x)^2$ and let
$X+X_0$ be represented by a quadruple
$\boldhats=(\uhat_0,\uhat_1,\vhat_0,\vhat_1)$ corresponding to
polynomials $\hat U(x)=x^2+\uhat_1\,x+\uhat_0$ and $\hat
V(x)=\vhat_1\,x+\vhat_0$ and such that $\hat U(x) \vert f(x)-\hat
V(x)^2$. Then for $i=1\,\ldots,9$, the equations
$g_i(\bolds,\boldhats)=0$ are satisfied.
\endproclaim

\head
A Biprojective Variety
\endhead

In general, a mapping of $C$ into $\Pic$ can be constructed for any $Q\in C$,
not just $Q=P_\infty$.  Thus for any $Q\in C$, there is a mapping 
$\phi_Q$ that imbeds $C$ into $\Pic$ given by
$$\phi_Q:\, P \mapsto 
\cases
Cl(P-Q)& \hbox{if $P \ne Q$;}\cr 0& \hbox{if $P=Q$;}\cr 
\endcases$$
where $Cl(D)$ denotes the divisor class of the divisor $D$.
The image of $C$ under $\phi_Q$ is just the translation of $\Theta$
by $Cl(P_\infty-Q)$, and is denoted $\Theta_Q$.  In particular
$\Theta=\Theta_{P_\infty}$.  For notational convenience, it is useful
to write $\Theta_{X_0}$ in place of $\Theta_{P_0}$.

For points $Q_1,Q_2\in C$ with $Q_1\ne Q_2$, it is interesting to
compute the intersection $\Theta_{Q_1}\cap\Theta_{Q_2}$.  This
consists of those divisor classes $G\in\Pic$ that contain two divisors
$D_1,D_2\in G$ of the forms $D_1 = P_1 - Q_1$ and $D_2 = P_2 - Q_2$
for some pair of points $P_1,P_2\in C$.  Since $D_1$ and $D_2$ are in
the same class,  they are equivalent, $D_1 \sim D_2$,  and therefore
$P_1 - P_2 \sim Q_1 - Q_2$.  Since $P_2 + \iota(P_2) - 2\cdot P_\infty
\sim 0$ and $Q_2 + \iota(Q_2) - 2\cdot P_\infty \sim 0$, this gives
$P_1+\iota(P_2)-2\cdot P_\infty \sim Q_1+\iota(Q_2)-2\cdot P_\infty$
and the assumption that $Q_2 \ne Q_2$ shows that the divisor on the
right is not in $\Theta$.  Now both of these divisors are reduced,
which means that the points must be the same, i.e. that either
$P_1=Q_1$ and $P_2=Q_2$ or $P_1=\iota(Q_2)$ and $Q_1=\iota(P_2)$.
The first of these possibilities gives $D_1=D_2=0$, while the second
case gives $\iota(Q_2)-Q_1 = D_1 \sim D_2 = \iota(Q_1) - Q_2$.
Therefore, 
$$\Theta_{Q_1}\cap\Theta_{Q_2}=\left\{0,Cl(\iota(Q_1)-Q_2)\right\}$$
and there are just two points in $\Theta_{Q_1}\cap\Theta_{Q_2}$.

In particular, this can be applied to find
$$\Theta \cap \Theta_{X_0} = \left\{0,Cl(P_0-P_\infty)\right\}$$
since $P_0$ is a Weierstrass point and $P_0=\iota(P_0)$.  Note
incidentally, that $X_0=Cl(P_0-P_\infty)$ is a point of order 2 in $\Pic$.

Based on this, it is useful to have a slightly finer classification of
divisor classes than the one above.  Every divisor class
$X\in\Pic$ contains a unique reduced divisor $D_X$ of one of the following
five forms:
\vskip 2pt
(ia) $D_X = P_1+P_2-2\cdot P_\infty$ with $P_1,P_2 \in \Caff$ 
and $P_1 \ne i(P_2)$ and $P_1,P_2 \ne P_0$;
\vskip 2pt
(ib) $D_X = P + P_0 - 2\dot P_\infty$ with $P \in \Caff$ and $P\ne P_0$;
\vskip 2pt
(iia) $D_X = P - P_\infty$ with $P \in \Caff$ and $P\ne P_0$;
\vskip 2pt
(iib) $D_X = P_0 - P_\infty$;
\vskip 2pt
(iii) $D_X = 0$.
\vskip2pt\noindent
Now cases (ia) and (ib) combine to form case (i) above and cases (iia)
and (iib) combine to form cases (ii) above.  With this finer
classification, elements of $\Pic$ that are neither in $\Theta$ nor in
$\Theta_{X_0}$ have reduced divisors of type (ia), elements of $\Theta_0$
that are not in $\Theta$ have reduced divisors of type (ib), elements
of $\Theta$ that are not in $\Theta_{X_0}$ have reduced divisors of type
(iia), and the two elements of $\Pic$ that are in both $\Theta$ and
$\Theta_{X_0}$ have reduced divisors of types (iib) and (iii).

\subhead
Glue Equations in Biprojective Space
\endsubhead

The idea now is to consider the ``graph'' of the map between the
affine varieties that parametrize $\Jac-\Theta$ and $\Jac-\Theta_0$ and look at
this as a biprojective variety.  The hope here is that except for
two points corresponding to cases (iib) and (iii) above, every point
will be accounted for.  In particular, the two varieties just need to
be projectivized with different homogenizing variables and then the
``glue'' needs to be applied which consists of bihomogenizing the
equations that relate the hatted and the unhatted variables.

Thus, there are projective coordinates $\bolds=(u_0,u_1,v_0,v_1,z)$ that
satisfy the pair of equations gotten by homogenizing the affine
equations for $\Jac-\Theta$ with respect to $z$.  Thus, the equations
are:
$$
\eqalign{
0 = E_0(\boldS) &= - a_2\,u_0\,z^3   + a_4\,{{u_0}^2}\,z^3 + a_3\,u_0\,u_1\,z^2 -
2\,{{u_0}^2}\,u_1\,z - a_4\,u_0\,{{u_1}^2}\,z + u_0\,{{u_1}^3} \nl
- {{v_0}^2}\,z^2 + u_0\,{{v_1}^2}\,z \cr
0 = E_1(\boldS) &= z^4 - a_3\,u_0\,z^3 + {{u_0}^2}\,z^2 - a_2\,u_1\,z^3 +
2\,a_4\,u_0\,u_1\,z^2 + a_3\,{{u_1}^2}\,z^2 \nl - 3\,u_0\,{{u_1}^2}\,z 
- a_4\,{{u_1}^3}\,z + {{u_1}^4} - 2\,v_0\,v_1\,z^2 + u_1\,{{v_1}^2}\,z \cr
}
$$
and it won't be necessary to go to the full projective closure.  In a similar
fashion, there are projective coordinates 
$\boldhatS=(\hat u_0,\hat u_1,\hat v_0,\hat v_1,\hat z)$ for $\Jac-\Theta_0$.
In this case, the equations are now the same as above but with all the
variables ``hatted''.  Thus, the equations are:
$$
\eqalign{
0 = \hat E_0(\boldhatS) &= - a_4\,\hat u_0\,\hat z^3   + a_2\,{{\hat u_0}^2}\,\hat z^3 
+ a_3\,\hat u_0\,\hat u_1\,\hat z^2 - 2\,{{\hat u_0}^2}\,\hat u_1\,\hat z \nl
- a_2\,\hat u_0\,{{\hat u_1}^2}\,\hat z 
+ \hat u_0\,{{\hat u_1}^3} - {{\hat v_0}^2}\,\hat z^2 +
\hat u_0\,{{\hat v_1}^2}\,\hat z \cr
0 = \hat E_1(\boldhatS) &= \hat z^4 - a_3\,\hat u_0\,\hat z^3 + {{\hat u_0}^2}\,\hat z^2 
- a_4\,\hat u_1\,\hat z^3 + 2\,a_2\,\hat u_0\,u_1\,\hat z^2 +
a_3\,{{\hat u_1}^2}\,\hat z^2 - 3\,\hat u_0\,{{\hat u_1}^2}\,\hat z \nl
- a_2\,{{\hat u_1}^3}\,\hat z + {{\hat u_1}^4} -
2\,\hat v_0\,\hat v_1\,\hat z^2 + \hat u_1\,{{\hat v_1}^2}\,\hat z \cr
}
$$
and again, going to the full projective closure is unnecessary.  Note
that setting $z=1$ in the first pair of equations recovers the affine
variety that they came from and that setting $\hat z=1$ in the second 
pair of equations recovers the affine variety that they came from.

The glue equations become the bihomogenization of the equations that
relate the unhatted variables to the hatted ones.  These equations are
$$\eqalign{
0 = G_1(\boldS,\boldhatS) &= u_0\,\hat u_0 - a_1\,z\,\hat z \cr
0 = G_2(\boldS,\boldhatS) &= \hat v_0\,u_0 + v_0\,\hat u_0 \cr
0 = G_3(\boldS,\boldhatS) &= v_0\,\hat v_0 + a_1\left(a_4\,z\,\hat z  - u_1\,\hat z - \hat u_1\,z\right) \cr
0 = G_4(\boldS,\boldhatS) &= \hat u_0\left(v_1\,\hat z+\hat v_1\,z\right) + \hat v_0\left(u_1\,\hat z - \hat u_1\,z\right) \cr
0 = G_5(\boldS,\boldhatS) &= u_0\left(v_1\,\hat z+\hat v_1\,z\right) + v_0\left(\hat u_1\,z - u_1\,\hat z\right) \cr
0 = G_6(\boldS,\boldhatS) &= a_1\,\hat u_0\,z^2 - a_1\,a_3\,z^2\,\hat z + a_1\,u_0\,\hat z\,z + a_1\,\hat u_1\,u_1\,z -
2\,u_1\,\hat v_0\,v_0 - 2\,\hat u_0\,v_0\,v_1 \cr
0 = G_7(\boldS,\boldhatS) &= a_1\,\hat u_0\,z\,\hat z - a_1\,a_3\,z\,\hat z^2 + a_1\,u_0\,\hat z^2 + a_1\,\hat u_1\,u_1\,\hat z -
2\,\hat u_1\,\hat v_0\,v_0 - 2\,u_0\,\hat v_0\,\hat v_1 \cr
0 = G_8(\boldS,\boldhatS) &= -a_1\,\hat v_0\,z^3 + a_3\,u_0\,\hat v_0\,z^2 - u_0^2\,\hat v_0\,z -
2\,a_4\,u_0\,u_1\,\hat v_0\,z + 3\,u_0\,u_1^2\,\hat v_0 
\cr&\quad + a_1\,u_1\,\hat v_1\,z^2 +
a_1\,u_1\,v_1\,z\,\hat z + 2\,\hat v_0\,v_0\,v_1\,z \cr
0 = G_9(\boldS,\boldhatS) &= -a_1\,v_0\,\hat z^3 + a_3\,\hat u_0\,v_0\hat z^2 - \hat u_0^2\,v_0\hat z -
2\,a_4\,\hat u_0\,\hat u_1\,v_0\,\hat z + 3\,\hat u_0\,\hat u_1^2\,v_0 
\cr&\quad +a_1\,\hat u_1\,v_1\,\hat z^2 + a_1\,\hat u_1\,\hat v_1\,\hat z\,z + 2\,v_0\,\hat v_0\,\hat v_1\,\hat z \cr
}$$
and it is worth noting that $G_1$, $G_2$, and $G_3$ are symmetric,
i.e. that $G_1(\boldS,\boldhatS)=G_1(\boldhatS,\bolds)$,
$G_2(\boldS,\boldhatS)=G_2(\boldhatS,\boldS)$, and
$G_3(\boldS,\boldhatS)=G_3(\boldhatS,\boldS)$, while $(G_4,G_5)$,
$(G_6,G_7)$, and $(G_8,G_9)$ are complementary pairs, i.e. 
$G_5(\boldS,\boldhatS)=G_4(\boldhatS,\boldS)$, 
$G_7(\boldS,\boldhatS)=G_6(\boldhatS,\boldS)$, and
$G_9(\boldS,\boldhatS)=G_8(\boldhatS,\boldS)$. It is also worth
noting that the homogeneous bidegrees of $G_1$, $G_2$ and $G_3$ are $(1,1)$,
the homogeneous bidegrees of $G_4$ and $G_7$ are $(1,2)$, the homogeneous bidegrees of
$G_5$ and $G_6$ are $(2,1)$, and homogeneous the bidegrees of $G_8$ and $G_9$
are $(3,1)$ and $(1,3)$, respectively, all of which follow from the inhomogeneous case.

Now consider what happens at the different types of divisor classes
and how they correspond to $z$ and $\hat z$ being nonzero or zero,
there being four cases.  The idea is that they should correspond to
the refined classifation of divisor classes.

If neither $z$ nor $\hat z$ is zero, then any biprojective point on
this whole set of equations corresponds to a unique point in $\Pic$ of
type (ia), and every divisor class of type (ia) gives rise to a unique
point on this biprojective variety with $z,\hat z \not= 0$.

Things get more interesting if either $z=0$ or $\hat z=0$, but not
both.  This should correspond to divisor classes of type (ib) and
(iia), and in fact this will be seen from looking at the equations,
but for the moment this correspondence is not assumed.

\subhead
The case $(z,\zhat)=(1,1)$
\endsubhead

In this case, the biprojective equations simply become the affine equations
that relate the coordinates of $X$ to the coordinates of $X+X_0$, as long
as neither $X$ nor $X+X_0$ are on $\Theta$. Alternatively, these can be thought 
of as equations on the affine variety $(\Jac-\Theta)\times(\Jac-\Theta_{X_0})$
where $\Theta_{X_0}$ is just the translation of $\Theta$ by the point $X_0$. Geometrically,
$\Jac-\Theta \simeq \Jac-\Theta_{X_0}$.

\subhead
The case $(z,\zhat)=(1,0)$
\endsubhead

Points of infinity type $(z,\zhat)=(0,1)$ should correspond to
divisors of type $P+P_0-2\,P_\infty$ for some affine $P\in\Caff$ with
$P \ne P_0$.  Similarly, points of infinity type $(z,\zhat)=(0,1)$,
should correspond to divisors of type $P-P_\infty$ for some affine
$P\in\Caff$ with $P \ne P_0$.

So suppose $\hat z=0$ and $z \not= 0$, in which case there is no loss
of generality in taking $z=1$.  Then it follows that $\hat u_1=0$ as
well.  The connecting equations after setting $z=1$ and $\hat z = 0$
then become (in order)
$$\eqalign{
G_1:\,\, 0 &= u_0\,\hat u_0 \cr
G_2:\,\, 0 &= \hat v_0\,u_0 + v_0\,\hat u_0 \cr
G_3:\,\, 0 &= v_0\,\hat v_0 + a_1\,\hat u_1 \cr
G_4:\,\, 0 &= \hat u_0\,\hat v_1 - \hat v_0\,\hat u_1 \cr
G_5:\,\, 0 &= u_0\,\hat v_1 + v_0\,\hat u_1 \cr
G_6:\,\, 0 &= a_1\,\hat u_0 + a_1\,\hat u_1\,u_1 - 2\,u_1\,\hat v_0\,v_0 - 2\,\hat u_0\,v_0\,v_1 \cr
G_7:\,\, 0 &= -2\,\hat u_1\,\hat v_0\,v_0 - 2\,u_0\,\hat v_0\,\hat v_1 \cr
G_8:\,\, 0 &= a_3\,u_0\,\hat v_0 - a_1\,\hat v_0 - u_0^2\,\hat v_0 -
              2\,a_4\,u_0\,u_1\,\hat v_0 + 3\,u_0\,u_1^2\,\hat v_0 + a_1\,u_1\,\hat v_1 + 2\,\hat v_0\,v_0\,v_1 \cr
G_9:\,\, 0 &= 3\,\hat u_0\,\hat u_1^2\,v_0\cr
}$$
from which some interesting implications will follow. 
Setting $\hat z = 0$ in equations $\hat E_0$ and $\hat E_1$ gives
$$
\eqalign{
E_0:\,\, 0 &= \hat u_0\,{{\hat u_1}^3} \cr
E_1:\,\, 0 &= {{\hat u_1}^4} \cr
}
$$
from which it follows that $\hat z = 0$ implies $\hat u_1 = 0$.
Upon setting $\hat u_1=0$, the glue equations now become (in order)
$$\eqalign{
G_1:\,\, 0 &= u_0\,\hat u_0 \cr
G_2:\,\, 0 &= \hat v_0\,u_0 + v_0\,\hat u_0 \cr
G_3:\,\, 0 &= v_0\,\hat v_0 \cr
G_4:\,\, 0 &= \hat u_0\,\hat v_1 \cr
G_5:\,\, 0 &= u_0\,\hat v_1 \cr
G_6:\,\, 0 &= a_1\,\hat u_0 - 2\,u_1\,\hat v_0\,v_0 - 2\,\hat u_0\,v_0\,v_1 \cr
G_7:\,\, 0 &= - 2\,u_0\,\hat v_0\,\hat v_1 \cr
G_8:\,\, 0 &= a_3\,u_0\,\hat v_0 - a_1\,\hat v_0 - u_0^2\,\hat v_0 -
     2\,a_4\,u_0\,u_1\,\hat v_0 + 3\,u_0\,u_1^2\,\hat v_0 + a_1\,u_1\,\hat v_1 + 2\,\hat v_0\,v_0\,v_1 \cr
}$$
Now suppose that $u_0 \not= 0$.  Then the first of these equations gives
$\hat u_0 = 0$ and the second equation becomes $0=\hat v_0\,u_0$ from
which it follows that $\hat v_0=0$.  The fifth equation now gives
$\hat v_1 = 0$.  So the assumption $u_0 \not= 0$ requires all the hatted
variables to be 0, but this is not a point in projective space, and
so does not correspond to a point on the biprojective variety.
Therefore $u_0=0$.  Now look at the above glue equations and note that
(in order) they become
$$\eqalign{
G_2:\,\, 0 &= v_0\,\hat u_0 \cr
G_3:\,\, 0 &= v_0\,\hat v_0 \cr
G_4:\,\, 0 &= \hat u_0\,\hat v_1 \cr
G_6:\,\, 0 &= a_1\,\hat u_0 - 2\,u_1\,\hat v_0\,v_0 - 2\,\hat u_0\,v_0\,v_1 \cr
G_8:\,\, 0 &= -a_1\,\hat v_0 + a_1\,u_1\,\hat v_1 + 2\,\hat v_0\,v_0\,v_1 \cr
}$$
and setting $z = 1$ and $u_0 = 0$ in equations $E_0$ and $E_1$ gives
$$
\eqalign{
E_0:\,\, 0 &= - {{v_0}^2} \cr
E_1:\,\, 0 &=  1 - a_2\,u_1 + a_3\,{{u_1}^2} - a_4\,{{u_1}^3} + {{u_1}^4} - 2\,v_0\,v_1 + u_1\,{{v_1}^2} \cr
}$$
so that $v_0=0$ from the first of these and then the second equation
now becomes
$$E_1:\,\, 0 =  1 - a_2\,u_1 + a_3\,{{u_1}^2} - a_4\,{{u_1}^3} + {{u_1}^4} + u_1\,{{v_1}^2}$$
and note that this implies $u_1 \not= 0$.  It is interesting to
multiply by $-u_1$ and rewrite this as
$$\left(u_1\,v_1\right)^2 = \left(-u_1\right)^5 +
a_4\,\left(-u_1\right)^4 + a_3\,\left(-u_1\right)^3 
+ a_2\,\left(-u_1\right)^2 + \left(-u_1\right)$$
so that the pair $\left(-u_1,u_1\,v_1\right)$ actually defines a point
on the punctured affine curve $\Caff-P_0$.

Also setting $v_0=0$ in the glue equations gives
$$\eqalign{
G_4:\,\, 0 &= \hat u_0\,\hat v_1 \cr
G_6:\,\, 0 &= a_1\,\hat u_0 \cr
G_8:\,\, 0 &= -a_1\,\hat v_0 + a_1\,u_1\,\hat v_1 \cr
}$$
and since $a_1 \not= 0$ it follows that $\hat u_0=0$ and the only
remaining glue equation is now
$$0 = -\hat v_0 + u_1\,\hat v_1$$ 
and since $u_1 \not= 0$ this gives a single point in the
projective space defined by  the hatted variables.  Thus
$${\hat v_0 \over \hat v_1} = u_1$$
so the projectived point on the hatted side of things is just
(minus) the $x$ coordinate of an affine point on $C$ (excluding
$P_0$).  This corresponds precisely to divisor classes of type
(ib), i.e. to classes having canonical representatives of the type
$P+P_0-2\cdot P_\infty$ where $P\in\Caff$ with $p \not= P_0$.

\proclaim{Proposition}
Suppose that $0$ is a root of $f(x)$ so $P_0=(0,0)$ is an affine
Weierstrass point on $C$. Let ${\boldkey P} =
(u_0,u_1,v_0,v_1,z)$ be the coordinates of a point on
$\hat A(\boldS)$, the projective closure of $A(\bolds)$, and let
${\hat\boldkey P} = (\uhat_0,\uhat_1,\vhat_0,\vhat_1,\zhat)$ be
the coordinates of a point on $\hat A(\hat\boldS)$, the projective
closure of $A(\hat\bolds^\prime)$. Furthermore also assume that these
coordinates also satisfy the bihomogenous glue equations
$\boldG(\boldS,\boldhatS,\bolda)=\boldzero$.  Suppose that $z=1$ and
$\zhat=0$. Then $u_0=v_0=\uhat_0=\uhat_1$ as well
and $\vhat_0=u_1\,\vhat_1$.  Furthermore $(-u_1,u_1v_1)$ is an
affine point on $C$ with $u_1\ne0$.
\endproclaim

Thus, points on the biprojective variety with infinity type
$(z,\zhat)=(1,0)$ are in one-to-one correpondence with divisor classes
of type (ib), as expected.

In a completely similar fashion, taking $z=0$ and $\hat z \not= 0$
gives divisor classes of type (iia).  In particular, in this case it
will be found that $u_0=\hat u_0=\hat v_0=u_1=0$ and that $\left(-\hat
u_1,\hat u_1\,\hat v_1\right)$ are the coordinates of a point on the
punctured affine curve $\Caff-P_0$ (so that $u_1 \not= 0$ and that on
the unhatted side of things, the remaining nonzero variables $v_0$ and
$v_1$ are in the fixed ratio
$${v_0 \over v_1} = \hat u_1.$$

\subhead
The case $(z,\zhat)=(0,0)$
\endsubhead

There are only two types of divisor classes left, namely (iib) and
(iii), which each consist of a single point, with representatives
$P_0-P_\infty$ and $0$, respectively.  These should correspond to 
the case $z=\hat z=0$.  Unfortunately, here is where there are some
problems.  Setting $z=\hat z=0$ in equations $E_1$ and $\hat E_1$
gives
$$
\eqalign{
E_1:\,\, 0 &= {u_1}^4 \cr
\hat E_1:\,\, 0 &= {\hat u_1}^4 \cr
}
$$
so that here $u_1=\hat u_1=0$.  Now setting $z=\hat z=u_1=\hat u_1=0$
in the glue equations gives (in order)
$$\eqalign{
G_1:\,\, 0 &= u_0\,\hat u_0 \cr
G_2:\,\, 0 &= \hat v_0\,u_0 + v_0\,\hat u_0 \cr
G_3:\,\, 0 &= v_0\,\hat v_0 \cr
G_6:\,\, 0 &= - 2\,\hat u_0\,v_0\,v_1 \cr
G_7:\,\, 0 &= - 2\,u_0\,\hat v_0\,\hat v_1 \cr
}$$
Also going through the same process of closure as in the affine case would give
$$0 = {{u_0}^2} {{v_1}^2}$$
and
$$0 = {{\hat u_0}^2} {{\hat v_1}^2}$$
which could also be of use if needed.

Now look what happens if $u_0\not =0$.  Then $\hat u_0=0$ from the
first equation and then $\hat v_0=0$ from the second equation.
Therefore in order that not all the hatted variables be zero this
gives $\hat v_1 \not= 0$, so the hatted variables all define a single
projective point.  For the unhatted variables $v_0$ is unconstrained,
but $v_1$ is forced to be 0. Therefore there is a projective line in
the unhatted variables. Similarly if $\hat u_0 \not =0$, the unhatted
variables all reduce to a single point, but the hatted variables give
a projective line.

If $u_0=\hat u_0 = 0$, then the only equation left is $0 = v_0\,\hat
v_0$ so either $v_0=0$ or $ \hat v_0=0$.  If $v_0=0$ then $v_1 \not=0$
or else all the unhatted variables would be 0, so that the unhatted
variables give a single projective point, and the hatted variables
$\hat v_0$ and $\hat v_1$ are unconstrained, so there is a projective
line in the hatted varaibles.  Similarly, if $\hat v_0=0$, then there
is a projective point in the hatted variables and a projective line in
the unhatted variables.

\proclaim{Proposition}
Suppose that $0$ is a root of $f(x)$ so $P_0=(0,0)$ is an affine
Weierstrass point on $C$. Let ${\boldkey P} =
(u_0,u_1,v_0,v_1,z)$ be the coordinates of a point on
$\hat A(\boldS)$, the projective closure of $A(\bolds)$, and let
${\hat\boldkey P} = (\uhat_0,\uhat_1,\vhat_0,\vhat_1,\zhat)$ be
the coordinates of a point on $\hat A(\hat\boldS)$, the projective
closure of $A(\hat\bolds^\prime)$. Furthermore also assume that these
coordinates also satisfy the bihomogenous glue equations
$\boldG(\boldS,\boldhatS,\bolda)=\boldzero$.  Suppose that $z=0$ and
$\zhat=0$. Then $u_1=\uhat_1=0$ and exactly one of the following occurs:
\vskip 1pt
(i) $v_1=\uhat_0=\vhat_0=0$ and $\vhat_1\ne0$ and $(u_0,v_0)\ne(0,0)$;
\vskip 1pt
(ii) $\vhat_1=u_0=v_0=0$ and $v_1\ne0$ and $(\uhat_0,\vhat_0)\ne(0,0)$;
\vskip 1pt
(iii) $u_0=\uhat_0=v_0=0$ and $v_1\ne0$ and $(\vhat_0,\vhat_1)\ne(0,0)$;
\vskip 1pt
(iv) $u_0=\uhat_0=\vhat_0=0$ and $\vhat_1\ne0$ and $(v_0,v_1)\ne(0,0)$;
\vskip 1pt\noindent
Each of these cases determines a biprojective line.

\endproclaim

Thus the case $z=\hat z=0$ gives a union of four projective lines
rather that a pair of projective points.  
The unfortunate result of all of this is that for the infinity type
$(z,\zhat)=(0,0)$, what results is not just a pair of points. 
The problem is now to somehow fix this up.

\subhead
More General Weierstrass Points: Adding a 2-Division Point in General
\endsubhead

It seems worthwhile to slightly rewrite the formulas so far
in a bit more generality to allow for an affine branch point 
on $C$ other than $(0,0)$. Eventually, it will be important to
allow multiple affine branch points and these formulas will be necessary.
The actual modifications to the formulas are really fairly simple.
Assume that $\rtzero$ is a root of $f(x)$, the monic quintic polynomial,
and write the affine equation for the curve in the form
$$y^2 = f(x)
= (x-\rtzero)^5 + \azero^\prime_4(x-\rtzero)^4 + \azero^\prime_3(x-\rtzero)^3 
  + \azero^\prime_2(x-\rtzero)^2 + \azero^\prime_1(x-\rtzero) $$
and the affine branch point in question is $P_\rtzero=(\rtzero,0)$.

The $U$ and $V$ polynomials now can be written as 
$$\eqalign{
U(x;X) &= (x-\rtzero)^2+u^\prime_1\,(x-\rtzero)+u^\prime_0 \cr
V(x;X) &= v^\prime_1\,(x-\rtzero)+v^\prime_0 \cr
}$$
and the condition to be satisfied is still
$$U(x;X) \vert f(x) - V(x;X)^2$$
which can be fulfilled by writing
$$f(x) - V(x;X)^2 \equiv e^\prime_1\,(x-\rtzero)+e^\prime_0 \bmod U(x;X)$$
and then insisting that
$$\eqalign{
0 &= e^\prime_0 \cr&=   
-\azero^\prime_2 {u^\prime_0} + \azero^\prime_4 {u^\prime_0}^2 
+ \azero^\prime_3 {u^\prime_0} {u^\prime_1} - 2 {u^\prime_0}^2 {u^\prime_1} 
- \azero^\prime_4 {u^\prime_0} {u^\prime_1}^2 + {u^\prime_0} {u^\prime_1}^3 
- {v^\prime_0}^2 + {u^\prime_0} {v^\prime_1}^2
\cr
0 &= e^\prime_1 \cr&= 
\azero^\prime_1 - \azero^\prime_3 {u^\prime_0} + {u^\prime_0}^2 
- \azero^\prime_2 {u^\prime_1} + 2 \azero^\prime_4 {u^\prime_0} {u^\prime_1} 
+ \azero^\prime_3 {u^\prime_1}^2 - 3 {u^\prime_0} {u^\prime_1}^2 
- \azero^\prime_4 {u^\prime_1}^3 + {u^\prime_1}^4 
\cr&\quad - 2 {v^\prime_0} {v^\prime_1} + {u^\prime_1} {v^\prime_1}^2
\cr
}$$
so that in these new coordinates, any quadruple
$\left({u^\prime_1},{u^\prime_0},{v^\prime_1},{v^\prime_0}\right)$ for
which both equations
$e^\prime_0\left({u^\prime_1},{u^\prime_0},{v^\prime_1},{v^\prime_0}\right)
= 0$ and 
$e^\prime_1\left({u^\prime_1},{u^\prime_0},{v^\prime_1},{v^\prime_0}\right)
= 0$ hold gives rise to a pair of polynomials $U(x) =
(x-\rtzero)^2+{u^\prime_1}\,(x-\rtzero)+{u^\prime_0}$ and $V(x) =
{v^\prime_1}\,(x-\rtzero)+{v^\prime_0}$ such that $U(x)\vert f(x)-V(x)^2$.
These two modified equations then define the same affine variety as before, 
$\Pic-\Theta$.  

Of course it is quite straightforward to convert
these equations to the old coordinates.  Just write
$$\eqalign{
U(x;X) &= (x-\rtzero)^2+{u^\prime_1}\,(x-\rtzero)+{u^\prime_0} \cr
&= x^2 + ({u^\prime_1} - 2\,\rtzero)\,x 
   + ({u^\prime_0} - {u^\prime_1}\,\rtzero + \rtzero^2) \cr
&= x^2+u_1\,x+u_0 \cr
V(x;X) &= {v^\prime_1}\,(x-\rtzero)+{v^\prime_0} \cr
&= {v^\prime_1}\,x+({v^\prime_0}-{v^\prime_1}\,\rtzero)\cr
&= v_1\,x+v_0 \cr
}$$
so that
$$\eqalign{
{u^\prime_1} &= u_1 + 2\,\rtzero \cr
{u^\prime_0} &= u_0 + u_1\,\rtzero + \rtzero^2 \cr
{v^\prime_1} &= v_1 \cr
{v^\prime_0} &= v_0 + v_1\,\rtzero \cr
}$$
which can be back substituted into equations $e^\prime_0=0$ and
$e^\prime_1=0$. 

It should be noted that with the substitutions
$$\eqalign{f(x)
&= x^5+a_4\,x^4+a_3\,x^3+a_2\,x^2+a_1\,x+a_0 \cr
&= (x-\rtzero)^5 + \azero^\prime_4(x-\rtzero)^4 + \azero^\prime_3(x-\rtzero)^3 
  + \azero^\prime_2(x-\rtzero)^2 + \azero^\prime_1(x-\rtzero) + \azero^\prime_0 \cr
}$$
it follows that
$$\eqalign{
a_4 &= \azero^\prime_4 - 5\,\rtzero \cr
a_3 &= \azero^\prime_3 - 4\,\azero^\prime_4\,\rtzero + 10\,\rtzero^2 \cr
a_2 &= \azero^\prime_2 - 3\,\azero^\prime_3\,\rtzero + 6\,\azero^\prime_4\,\rtzero^2 
       - 10\,\rtzero^3 \cr
a_1 &= \azero^\prime_1 - 2\,\azero^\prime_2\,\rtzero + 3\,\azero^\prime_3\,\rtzero^2 
       - 4\,\azero^\prime_4\,\rtzero^3 + 5\,\rtzero^4 \cr
a_0 &= \azero^\prime_0 - \azero^\prime_1\,\rtzero + \azero^\prime_2\,\rtzero^2 
       - \azero^\prime_3\,\rtzero^3 + \azero^\prime_4\,\rtzero^4 - \rtzero^5 \cr
}$$
which is a substitution that must be made when using Weierstrass points other than $(0,0)$
(in which case $\azero^\prime_0=0$). In matrix form this is
$$
\pmatrix
a_0 \cr a_1 \cr a_2 \cr a_3 \cr a_4 \cr
\endpmatrix
=
\pmatrix
1 & -\rtzero & \rtzero^2    & -\rtzero^3   & \rtzero^4     \cr
0 & 1        & -2\,\rtzero  & 3\,\rtzero^2 & -4\,\rtzero^3 \cr
0 & 0        & 1            & -3\,\rtzero  & 6\,\rtzero^2  \cr
0 & 0        & 0            & 1            & -4\,\rtzero   \cr
0 & 0        & 0            & 0            & 1             \cr
\endpmatrix
\,
\pmatrix
\azero^\prime_0 \cr \azero^\prime_1 \cr \azero^\prime_2 \cr \azero^\prime_3 \cr \azero^\prime_4 \cr
\endpmatrix
+
\pmatrix
- \rtzero^5 \cr 5\,\rtzero^4 \cr - 10\,\rtzero^3 \cr 10\,\rtzero^2 \cr - 5\,\rtzero \cr
\endpmatrix
$$
and the inverse transformation is 
$$
\pmatrix
\azero^\prime_0 \cr \azero^\prime_1 \cr \azero^\prime_2 \cr \azero^\prime_3 \cr \azero^\prime_4 \cr
\endpmatrix
=
\pmatrix
1 & \rtzero & \rtzero^2   & \rtzero^3    & \rtzero^4    \cr
0 & 1       & 2\,\rtzero  & 3\,\rtzero^2 & 4\,\rtzero^3 \cr
0 & 0       & 1           & \rtzero      & 6\,\rtzero^2 \cr
0 & 0       & 0           & 1            & 4\,\rtzero   \cr
0 & 0       & 0           & 0            & 1            \cr
\endpmatrix
\,
\pmatrix
a_0 \cr a_1 \cr a_2 \cr a_3 \cr a_4 \cr
\endpmatrix
+
\pmatrix
\rtzero^5 \cr 5\,\rtzero^4 \cr 10\,\rtzero^3 \cr 10\,\rtzero^2 \cr 5\,\rtzero \cr
\endpmatrix
.$$

It is also interesting to note that the relationship between the defining polynomials
$(e^\prime_1,e^\prime_0)$ and $(e_1,e_0)$ is also given by 
$$\eqalign{
e^\prime_1\left({u^\prime_1},{u^\prime_0},{v^\prime_1},{v^\prime_0}\right)
&=e_1(u_1,u_0,v_1,v_0) \cr
e^\prime_0\left({u^\prime_1},{u^\prime_0},{v^\prime_1},{v^\prime_0}\right)
&=e_0(u_1,u_0,v_1,v_0) + \rtzero\,e_1(u_1,u_0,v_1,v_0) \cr
}$$
which follows directly from substituting $x\mapsto x+\rtzero$ in
$$f(x)-V(x)^2 \equiv e_1\,x+e_0 \bmod U(x)$$ so that
$$f(x+\rtzero)-V(x+\rtzero) \equiv e_1\,(x+\rtzero)+e_0 
= e_1^\prime\,x+e_0^\prime \bmod U(x+\rtzero).$$

\head
A quadri-projective variety
\endhead

The next question is whether the problem at infinity with the
biprojective variety (which was obtained by working with a single
affine Weierstrass point)(which was obtained by working with a single
affine Weierstrass point) can be cleared up by working with a second
affine Weierstrass point. 

To that end, let $\rho^{(1)}$ and $\rho^{(2)}$ be distinct roots of
the quintic $f(x)$, so that $P_1=\bigl(\rho^{(1)},0\bigr)$ and
$P_2=\bigl(\rho^{(2)},0\bigr)$ are distinct Weierstrass points on the
projective curve $C$ defined by $y^2=f(x)$. The divisors
$P_1-P_\infty$ and $P_2-P_\infty$ determine divisor classes $X_1$ and
$X_2$ in $\Pic$.  Furthermore, the divisor of the function
$x-\rho^{(1)}$ is $2P_1-2P_\infty$, so the divisor class $X_1$ is of
order 2 in $\Pic$, and similarly the divisor class $X_2$ is also of
order 2.

With these two distinguished affine points, there is an even finer
classification of divisor classes than before. The new classification
is as follows: Every divisor class
$X\in\Pic$ contains a unique reduced divisor $D_X$ of one of the following
eight forms:
\vskip 2pt
(1a) $D_X = Q_1+Q_2-2\cdot P_\infty$ with $Q_1,Q2 \in \Caff$ 
and $Q_1 \ne i(Q_2)$ and $Q_1,Q_2 \notin \{P_1,P_2\}$;
\vskip 2pt
(1b) $D_X = Q+P_1-2\cdot P_\infty$ with $Q \in \Caff$ and
$Q \notin \{P_1,P_2\}$;
\vskip 2pt
(1c) $D_X = Q+P_2-2\cdot P_\infty$ with $Q \in \Caff$ and
$Q \notin \{P_1,P_2\}$;
\vskip 2pt
(1d) $D_X = P_1+P_2-2\cdot P_\infty$;
\vskip 2pt
(2a) $D_X = Q-P_\infty$ with $Q \in \Caff$ and
$Q \notin \{P_1,P_2\}$;
\vskip 2pt
(2b) $D_X = P_1-P_\infty$;
\vskip 2pt
(2c) $D_X = P_2-P_\infty$;
\vskip 2pt
(3) $0$.
\vskip2pt\noindent
Type (1a) is by far the largest of these types and comprises a two-dimensional
affine variety. Types (1b), (1c), and (2a) each correspond to one-dimensional
affine varieties that are isomorphic is $\Caff-\{P_1,P_2\}$. Types (1d), (2b), (2c),
and (3) each correspond to a single point.

Every divisor class is one of these types, and these types don't
overlap.  Furthermore, as will be seen in the following analysis, each
of these types of divisor class corresponds to a distinct infinity
type in the quadri-homogeneous variety constructed in this session,
with the exception of type (1a) which includes two different infinity
types. What needs to be checked is that there are no other infinity
types that contain actual geometric points.  Unfortunately, this will
not turn out to be the case. However, where this goes awry is with
the infinity type $(0,0,0,0)$, and the way that it goes wrong is much
less dramatic than the case of the bihomogeneous variety described above.

\subhead
Multi-homogeneous Coordinates
\endsubhead

In an effort to fix up the problems at infinity that occurred in the bihomogeneous
analysis, the next idea is to consider two rational affine Weierstrass points on the
genus 2 curve $C$. Thus $C$ is still defined by $y^2 = f(x)$ where $f(x)$ is a 
monic polynomial of degree 5 with no multiple roots.  Now suppose that $\rtone$
and $\rttwo$ are two distinct roots of $f(x)$, i.e. $f(\rtone)=f(\rttwo)=0$ with
$\rtone\ne\rttwo$ and in general write
$$\eqalign{f(x)
&=\bigl(x-\rho^{(1)}\bigr)\,\bigl(x-\rho^{(2)}\bigr)\,\bigl(x-\rho^{(3)}\bigr)
  \,\bigl(x-\rho^{(4)}\bigr)\,\bigl(x-\rho^{(5)}\bigr)\cr
&=\bolda_0\cdot\boldx=\bolda_{\rho^{(i)}}\cdot\boldx_{\rho^{(i)}}\cr
&=x^5+\azero_4\,x^4+\azero_3\,x^3+\azero_2\,x^2+\azero_1\,x+\azero_0 \cr
&=(x-\rho^{(i)})^5 \!+ \azero^{(i)}_4 (x-\rho^{(i)})^4 \!+ \azero^{(i)}_3 (x-\rho^{(i)})^3 
  \!+ \azero^{(i)}_2 (x-\rho^{(i)})^2 \!+ \azero^{(i)}_1 (x-\rho^{(i)}) \cr
}$$
for $i=1,\ldots,5$, where
$$\eqalign{
\bolda_0&=(\azero_0,\azero_1,\azero_2,\azero_3,\azero_4,1) \cr
\boldx&=(1,x,x^2,x^3,x^4,x^5) \cr
\bolda_{\rho^{(i)}}&=(0,\azero^{(i)}_1,\azero^{(i)}_2,\azero^{(i)}_3,\azero^{(i)}_4,1) \cr
\boldx_{\rho^{(i)}}&=(1,(x-\rho^{(i)}),(x-\rho^{(i)})^2,(x-\rho^{(i)})^3,(x-\rho^{(i)})^4,(x-\rho^{(i)})^5)  \cr
}$$
and let 
$$\eqalign{
\bolds&=(u_0,u_1,v_0,v_1)\cr
\boldhats&=(\uhat_0,\uhat_1,\vhat_0,\vhat_1)\cr
}$$
be sets of variables, and let
$$\eqalign{
\boldS&=(u_0,u_1,v_0,v_1,z)\cr
\boldhatS&=(\uhat_0,\uhat_1,\vhat_0,\vhat_1,\zhat)\cr
}$$
be the corresponding sets of homogeneous variables.

Now let
$$\eqalign{
&e_0(v_1,v_0,u_1,u_0) \cr 
&= \azero_0 - \azero_2\,u_0 + \azero_4\,u_0^2 + \azero_3\,u_0\,u_1 - 2\,u_0^2\,u_1 - \azero_4\,u_0\,u_1^2 + u_0\,u_1^3 - v_0^2 + u_0\,v_1^2
\cr
&e_1(v_1,v_0,u_1,u_0) \cr 
&= \azero_1 - \azero_3\,u_0 + u_0^2 - \azero_2\,u_1 + 2\,\azero_4\,u_0\,u_1 
   + \azero_3\,u_1^2 - 3\,u_0\,u_1^2 - \azero_ 4\,u_1^3 + u_1^4 \cr
   &\qquad - 2\,v_0\,v_1 + u_1\,v_1^2
\cr
}$$
so the affine variety defined by
$e_0(v_1,v_0,u_1,u_0)=e_1(v_1,v_0,u_1,u_0)=0$ is $\JacC-\Theta$.
It is useful to write
$\bolde(\bolds,\boldazero)=\bigl(e_0(\bolds,\boldazero),e_1(\bolds,\boldazero)\bigr)$
so the equations defining the affine variety are
$\bolde(\bolds,\boldazero)=\boldzero$.
The homogenized versions of these polynomials are
$$\eqalign{
&E_0(v_1,v_0,u_1,u_0,z) = E_0(\boldS,\boldazero) \cr 
&= \azero_0\,z^4 - \azero_2\,u_0\,z^3 + \azero_4\,u_0^2\,z^2 + \azero_3\,u_0\,u_1\,z^2 
 - 2\,u_0^2\,u_1\,z - \azero_4\,u_0\,u_1^2\,z + u_0\,u_1^3 \cr
 &\qquad- v_0^2\,z^2 + u_0\,v_1^2\,z \cr
&E_1(v_1,v_0,u_1,u_0,z) = E_1(\boldS,\boldazero) \cr 
&= \azero_1\,z^4 - \azero_3\,u_0\,z^3 + u_0^2\,z^2 - \azero_2\,u_1\,z^3 + 2\,\azero_4\,u_0\,u_1\,z^3 
 + \azero_3\,u_1^2\,z^2 - 3\,u_0\,u_1^2\,z \cr&\quad - \azero_ 4\,u_1^3\,z + u_1^4 
 - 2\,v_0\,v_1\,z^2 + u_1\,v_1^2\,z^2 \cr
}$$
and the projectve closure of the equations
$E_0(\boldS,\boldazero)=E_1(\boldS,\boldazero)=0$ is just the
projective closure of $\Jac-\Theta$.

In an effort to fix up the problems at infinity that occurred in the bihomogeneous
analysis, the next idea is to consider two rational affine Weierstrass points on the
genus 2 curve $C$. Thus $C$ is still defined by $y^2 = f(x)$ where $f(x)$ is a 
monic polynomial of degree 5 with no multiple roots.  Now suppose that $\rtone$
and $\rttwo$ are two distinct roots of $f(x)$, i.e. $f(\rtone)=f(\rttwo)=0$ with
$\rtone\ne\rttwo$ and in general write
$$\eqalign{f(x)
&=\bigl(x-\rho^{(1)}\bigr)\,\bigl(x-\rho^{(2)}\bigr)\,\bigl(x-\rho^{(3)}\bigr)
  \,\bigl(x-\rho^{(4)}\bigr)\,\bigl(x-\rho^{(5)}\bigr)\cr
&=\bolda_0\cdot\boldx=\bolda_{\rho^{(i)}}\cdot\boldx_{\rho^{(i)}}\cr
&=x^5+\azero_4\,x^4+\azero_3\,x^3+\azero_2\,x^2+\azero_1\,x+\azero_0 \cr
&=(x-\rho^{(i)})^5 \!+ \azero^{(i)}_4 (x-\rho^{(i)})^4 \!+ \azero^{(i)}_3 (x-\rho^{(i)})^3 
  \!+ \azero^{(i)}_2 (x-\rho^{(i)})^2 \!+ \azero^{(i)}_1 (x-\rho^{(i)}) \cr
}$$
for $i=1,\ldots,5$, where
$$\eqalign{
\bolda_0&=(\azero_0,\azero_1,\azero_2,\azero_3,\azero_4,1) \cr
\boldx&=(1,x,x^2,x^3,x^4,x^5) \cr
\bolda_{\rho^{(i)}}&=(0,\azero^{(i)}_1,\azero^{(i)}_2,\azero^{(i)}_3,\azero^{(i)}_4,1) \cr
\boldx_{\rho^{(i)}}&=(1,(x-\rho^{(i)}),(x-\rho^{(i)})^2,(x-\rho^{(i)})^3,(x-\rho^{(i)})^4,(x-\rho^{(i)})^5)  \cr
}$$
and let 
$$\eqalign{
\bolds&=(u_0,u_1,v_0,v_1)\cr
\boldhats&=(\uhat_0,\uhat_1,\vhat_0,\vhat_1)\cr
}$$
be sets of variables, and let
$$\eqalign{
\boldS&=(u_0,u_1,v_0,v_1,z)\cr
\boldhatS&=(\uhat_0,\uhat_1,\vhat_0,\vhat_1,\zhat)\cr
}$$
be the corresponding sets of homogeneous variables.

Now let
$$\eqalign{
&e_0(v_1,v_0,u_1,u_0) \cr 
&= \azero_0 - \azero_2\,u_0 + \azero_4\,u_0^2 + \azero_3\,u_0\,u_1 - 2\,u_0^2\,u_1 - \azero_4\,u_0\,u_1^2 + u_0\,u_1^3 - v_0^2 + u_0\,v_1^2
\cr
&e_1(v_1,v_0,u_1,u_0) \cr 
&= \azero_1 - \azero_3\,u_0 + u_0^2 - \azero_2\,u_1 + 2\,\azero_4\,u_0\,u_1 
   + \azero_3\,u_1^2 - 3\,u_0\,u_1^2 - \azero_ 4\,u_1^3 + u_1^4 \cr
   &\qquad - 2\,v_0\,v_1 + u_1\,v_1^2
\cr
}$$
so the affine variety defined by
$e_0(v_1,v_0,u_1,u_0)=e_1(v_1,v_0,u_1,u_0)=0$ is $\JacC-\Theta$.
It is useful to write
$\bolde(\bolds,\boldazero)=\bigl(e_0(\bolds,\boldazero),e_1(\bolds,\boldazero)\bigr)$
so the equations defining the affine variety are
$\bolde(\bolds,\boldazero)=\boldzero$.
The homogenized versions of these polynomials are
$$\eqalign{
&E_0(v_1,v_0,u_1,u_0,z) = E_0(\boldS,\boldazero) \cr 
&= \azero_0\,z^4 - \azero_2\,u_0\,z^3 + \azero_4\,u_0^2\,z^2 + \azero_3\,u_0\,u_1\,z^2 
 - 2\,u_0^2\,u_1\,z - \azero_4\,u_0\,u_1^2\,z + u_0\,u_1^3 \cr
 &\qquad- v_0^2\,z^2 + u_0\,v_1^2\,z \cr
&E_1(v_1,v_0,u_1,u_0,z) = E_1(\boldS,\boldazero) \cr 
&= \azero_1\,z^4 - \azero_3\,u_0\,z^3 + u_0^2\,z^2 - \azero_2\,u_1\,z^3 + 2\,\azero_4\,u_0\,u_1\,z^3 
 + \azero_3\,u_1^2\,z^2 - 3\,u_0\,u_1^2\,z \cr&\quad - \azero_ 4\,u_1^3\,z + u_1^4 
 - 2\,v_0\,v_1\,z^2 + u_1\,v_1^2\,z^2 \cr
}$$
and the projectve closure of the equations
$E_0(\boldS,\boldazero)=E_1(\boldS,\boldazero)=0$ is just the
projective closure of $\JacC-\Theta$, which (unfortunately) is not
$\JacC$. 

\subhead
Glue for Other Weierstrass Points
\endsubhead

Let 
$$M_{\rti}=
\pmatrix 1       & 0     & 0    & 0      & 0 \cr
        \rti     & 1     & 0    & 0      & 0 \cr
        0        & 0     & 1    & 0      & 0 \cr
        0        & 0     & \rti & 1      & 0 \cr
        {\rti}^2 & 2\,\rti & 0      & 0      & 1 \cr
\endpmatrix
$$
which effects the projective transformation
$$\eqalign{
{u}_0 &\mapsto u_0 + u_1\,\rti + z\,{\rti}^2 \cr
{u}_1 &\mapsto u_1 + 2\,z\,\rti \cr
{v}_0 &\mapsto v_0 + v_1\,\rti \cr
{v}_1 &\mapsto v_1 \cr
{z}   &\mapsto z \cr
}$$
when written as
$\boldS \mapsto \boldS\,M_{\rti}$. The sets of equations
$\boldE(\boldS\,M_{\rti},\boldaone)=\boldzero$ for $i=1,\ldots,5$
are equivalent to the set of equations
$\boldE(\boldS,\boldazero)=\boldzero$.

In terms of the inhomogeneous equations, let 
$$\eqalign{
m_{\rti}&=
\pmatrix 1     & 0    & 0    & 0      \cr
        \rti   & 1    & 0    & 0      \cr
        0      & 0    & 1    & 0      \cr
        0      & 0    & \rti & 1      \cr \endpmatrix \cr
n_{\rti}&=
\pmatrix {\rti}^2& 2\,\rti & 0      & 0\cr \endpmatrix \cr
}$$
so that $\bolds\mapsto \bolds\,m_{\rti}+n_{\rti}$ effects the 
corresponding inhomogeneous transformation
$$\eqalign{
{u}_0 &\mapsto u_0 + u_1\,\rti + {\rti}^2 \cr
{u}_1 &\mapsto u_1 + 2\,\rti \cr
{v}_0 &\mapsto v_0 + v_1\,\rti \cr
{v}_1 &\mapsto v_1 \cr
}$$
which will prove useful.

\subhead
Quadrihomogeneous Coordinates
\endsubhead

Let $P_{\rtone}=(\rtone,0)\in C$ and $P_{\rttwo}=(\rttwo,0)\in C$, and set
$X_{1}=\Cl(P_{\rtone}-P_\infty)\in\JacC$ and $X_{2}=\Cl(P_{\rttwo}-P_\infty)\in\JacC$,
so that $X_{1}$ and $X_{2}$ are distinct and nontrivial $[2]$-division points
in $\JacC$. If $\Theta_X$ denotes the translation of the $\Theta$ divisor on $\JacC$
by the point $X\in\JacC$ then there is a simple correspondence between affine varieties
$$\labeledsquare{\JacC-\Theta_{}}{\JacC-\Theta_{X_{1}}}{\JacC-\Theta_{X_{2}}}{\JacC-\Theta_{X_{1}+X_{2}}}{\scriptstyle(+X_{1})}{\scriptstyle(+X_{2})}{\scriptstyle(+X_{1})}{\scriptstyle(+X_{2})}
$$
and the reason that the arrows are bidirectional is that adding $X_{1}$ maps $\JacC-\Theta_{}$
to $\JacC-\Theta_{X_{1}}$ and also $\JacC-\Theta_{X_{1}}$ back to $\JacC-\Theta_{}$, and 
similarly for all the other arrows. The corresponding inhomogeneous equations are 
$$\labeledsquare{\bolde(\bolds_{\zz},\boldazero)}{\bolde(\bolds_{\zo},\boldazero)}{\bolde(\bolds_{\oz},\boldazero)}{\bolde(\bolds_{\oo},\boldazero)}{\scriptstyle\boldg(\bolds_{\zz}M_{\rtone},\bolds_{\zo}M_{\rtone},\bolda_{\rtone})}{\scriptstyle\boldg(\bolds_{\zz}M_{\rttwo},\bolds_2M_{\rttwo},\bolda_{\rttwo})}{\scriptstyle\boldg(\bolds_{\zo}M_{\rttwo},\bolds_{\oo}M_{\rttwo},\bolda_{\rttwo})}{\scriptstyle\boldg(\bolds_{\oz}M_{\rtone},\bolds_{\oo}M_{\rtone},\bolda_{\rtone})}$$
where
$$\bolds_{ij}=\bigl(u_{ij,0},u_{ij,1},v_{ij,0},v_{ij,1}\bigr)$$
for $ij\in\{\zz,\zo,\oz,\oo\}$. Here
$$\bolds_{ij}M_{\rtone}
=\bigl(u_{ij,0}+u_{ij,1}\rtone+{\rtone}^2,u_{ij,1}+2\rtone,v_{ij,0}+v_{ij,1}\rtone,v_{ij,1}\bigr)$$
with a similar formula for $\bolds_{ij}M_{\rttwo}$.

The corresponding homogeneous polynomials fit into the diagram
$$\labeledsquare{\boldE(\boldS_{\zz},\boldazero)}{\boldE(\boldS_{\zo},\boldazero)}{\boldE(\boldS_{\oz},\boldazero)}{\boldE(\boldS_{\oo},\boldazero)}{\scriptstyle\boldG(\boldS_{\zz}M_{\rtone},\boldS_{\zo}M_{\rtone},\bolda_{\rtone})}{\scriptstyle\boldG(\boldS_{\zz}M_{\rttwo},\boldS_2M_{\rttwo},\bolda_{\rttwo})}{\scriptstyle\boldG(\boldS_{\zo}M_{\rttwo},\boldS_{\oo}M_{\rttwo},\bolda_{\rttwo})}{\scriptstyle\boldG(\boldS_{\oz}M_{\rtone},\boldS_{\oo}M_{\rtone},\bolda_{\rtone})}
$$
where
$$\boldS_{ij}=\bigl(u_{ij,0},u_{ij,1},v_{ij,0},v_{ij,1},z_{ij}\bigr)$$
for $ij\in\{\zz,\zo,\oz,\oo\}$. Here
$$\boldS_{ij}M_{\rtone}
=\bigl(u_{ij,0}+u_{ij,1}\rtone+z_{ij}{\rtone}^2,u_{ij,1}+2z_{ij}\rtone,v_{ij,0}+v_{ij,1}\rtone,v_{ij,1},z_{ij}\bigr)$$
with a similar formula for $\boldS_{ij}M_{\rttwo}$.

\subhead
Quadrihomogeneous Coordinates: Analysis of Infinity Types
\endsubhead

In order to see what the multi-projective variety determined by these
equations is, it is really necessary to understand what happens at
infinity in each of these components given the glue that holds them
all together. This can be done most easily by starting with the case
$\rtone=0$ and then translating to arbitrary $\rtone$, and noting that
these same equations apply for $\rttwo$.

Let $ij$ and $i^\prime j^\prime$ be two different indices taken from the set $\{\zz,\zo,\oz,\oo\}$
Now with $z_{ij}=1$ and $z_{i^\prime j^\prime}=0$, (which also implies $=u_{i^\prime j^\prime,1}=0$),
$U_{ij}\mapsto\boldS_{ij}M_{\rtone}$ and $U_{i^\prime j^\prime}\mapsto\boldS_{i^\prime j^\prime}M_{\rtone}$
at $z_{ij}=1$, and $z_{i^\prime j^\prime}=u_{i^\prime j^\prime,1}=0$
effects the transformation
$$\eqalign{
u_0&\mapsto u_{ij,0}+{\rtone}\,u_{ij,1}+{\rtone}^2 \cr
u_1&\mapsto u_{ij,1}+2\,{\rtone} \cr
v_0&\mapsto v_{ij,0}+{\rtone}\,v_{ij,1} \cr
v_1&\mapsto v_{ij,1} \cr
\uhat_0 &\mapsto u_{i^\prime j^\prime,0}+{\rtone}\,u_{i^\prime j^\prime,1} \cr
\uhat_1 &\mapsto u_{i^\prime j^\prime,1} \cr
\vhat_0&\mapsto v_{i^\prime j^\prime,0}+{\rtone}\,v_{i^\prime j^\prime,1} \cr
\vhat_1&\mapsto v_{i^\prime j^\prime,1} \cr
}$$
and gives rise to the following relations for 
$\boldG(\boldS_{ij}M_{\rtone},\boldS_{i^\prime j^\prime}M_{\rtone},\bolda_{\rtone})$
$$\eqalign{
0&= u_{ij,0}+{\rtone}\,u_{ij,1}+{\rtone}^2\cr 
0&= v_{ij,0}+{\rtone}\,v_{ij,1} \cr
0&=u_{i^\prime j^\prime,0}+{\rtone}\,u_{i^\prime j^\prime,1} = u_{i^\prime j^\prime,0}\cr
0&=u_{i^\prime j^\prime,1} \cr
0&=-(v_{i^\prime j^\prime,0}+{\rtone}\,v_{i^\prime j^\prime,1}) + (u_{ij,1}+2\,{\rtone})\,v_{i^\prime j^\prime,1} \cr
 &=-v_{i^\prime j^\prime,0}+(u_{ij,1}+{\rtone})\,v_{i^\prime j^\prime,1} \cr
}$$
with similar relations where $\rtone$ is replaced by $\rttwo$ (when adding the
order 2 point $X_{2}$ instead of $X_{1}$).

Now in the variable sets $\boldS_{\zz}$ and $\boldS_{\zo}$ with
$z_{\zz}=z_{\zo}=0$, (which also imply $u_{\zz,1}=u_{\zo,1}=0$), the
three equations that result are
$$\eqalign{
0&=u_{\zz,0}\,u_{\zo,0} \cr
0&=-(v_{\zo,0}+{\rtone}v_{\zo,1})\,u_{\zz,0} + (v_{\zz,0}+{\rtone}v_{\zz,1})\,u_{\zo,0} \cr
0&=(v_{\zz,0}+{\rtone}v_{\zz,1})(v_{\zo,0}+{\rtone}v_{\zo,1}) \cr
}$$
which are the only relations that can be inferred from these glue equations.

The goal here is that every point of $\JacC$ corresponds to a unique
point in the multi-projective variety defined by these equations, and
conversely that every point of the multi-projective variety corresponds
to a unique point of $\JacC$. There are several cases to consider
depending on the different infinity types. There are a lot of
symmetries present which will reduce the number of cases that need to
be considered.

\subhead
The Case $\boldz=(1,1,1,1)$
\endsubhead

The simplest case is where all the homogenizing variables $z_{\zz}$,
$z_{\zo}$, $z_{\oz}$, $z_{\oo}$ are non-zero. This corresponds to divisor
classes $Z$ such that none of $Z$, $Z+X_{1}$, $Z+X_{2}$,
$Z+X_{1}+X_{2}$ are in $\Theta$, and it is therefore
possible to take $z_{\zz}=z_{\zo}=z_{\oz}=z_{\oo}=1$. This case presents no
problems at all. Pictorially, this is
$$\labeledsquare{z_{\zz}=1}{z_{\zo}=1}{z_{\oz}=1}{z_{\oo}=1}{\scriptstyle(+X_{1})}{\scriptstyle(+X_{2})}{\scriptstyle(+X_{1})}{\scriptstyle(+X_{2})}
$$
with all the the corners corresponding to affine points on (translated) copies of 
$\Jac(C)-\Theta$, the the glue equations are all affine.

\proclaim{Proposition}
Solutions to the system of equations $B_4$ with $\boldz=(1,1,1,1)$ are in
one-to-one correspondence with the set of division classes in 
$\Pic-\bigl(\Theta\cup\Theta_{X_{1}}\cup\Theta_{X_{2}}\cup\Theta_{X_{1}+X_{2}}\bigr)$.
\endproclaim

\subhead
The Case $\boldz=(0,1,1,1)$
\endsubhead

The next simplest case is where three of the homogenizing variables
are non-zero and the fourth is zero. This corresponds to divisor
classes $Z$ such that exactly one of $Z$, $Z+X_{1}$,
$Z+X_{2}$, $Z+X_{1}+X_{2}$ is in $\Theta$, the others
being in $\JacC-\Theta$. This also presents no problems. Pictorially, the four
cases are
$$\matrix
\labeledsquare{z_{\zz}=0}{z_{\zo}=1}{z_{\oz}=1}{z_{\oo}=1}{\scriptstyle(+X_{1})}{\scriptstyle(+X_{2})}{\scriptstyle(+X_{1})}{\scriptstyle(+X_{2})} &
\labeledsquare{z_{\zz}=1}{z_{\zo}=0}{z_{\oz}=1}{z_{\oo}=1}{\scriptstyle(+X_{1})}{\scriptstyle(+X_{2})}{\scriptstyle(+X_{1})}{\scriptstyle(+X_{2})} \cr
\vphantom{\Bigg\vert} \cr
\labeledsquare{z_{\zz}=1}{z_{\zo}=1}{z_{\oz}=0}{z_{\oo}=1}{\scriptstyle(+X_{1})}{\scriptstyle(+X_{2})}{\scriptstyle(+X_{1})}{\scriptstyle(+X_{2})} &
\labeledsquare{z_{\zz}=1}{z_{\zo}=1}{z_{\oz}=1}{z_{\oo}=0}{\scriptstyle(+X_{1})}{\scriptstyle(+X_{2})}{\scriptstyle(+X_{1})}{\scriptstyle(+X_{2})} \cr
\endmatrix
$$
and now by symmetry it is suffices to consider the case $z_{\zz}=0$
and $z_{\zo}=z_{\oz}=z_{\oo}=1$. By viewing each of the arrows as a
simple bihomogeneous case, and referring to the previous analysis, it
is readily apparent that the $z_{\zz}=0$ refers to an element of
$\Theta$, and in particular $u_{\zz,0}=v_{\zz,0}=0$ and
$\bigl(-u_{\zz,1},u_{\zz,1}v_{\zz,1}\bigr)$ are the coordinates of an
affine point $P$ on $C$ that is neither $P_1=(\rho^{(1)},0)$ nor
$P_2=(\rho^{(2)},0)$. Then the affine points on the corners
$z_{\zo}=1$, $z_{\oz}=1$, and $z_{\oo}=1$ represent the divisor
classes of $P+P_1-2\,P_\infty$, $P+P_2-2\,P_\infty$, and 
$P+P_1+P_2-3\,P_\infty$, respectively.

\proclaim{Proposition}
Solutions to the system of equations $B_4$ with $\boldz=(0,1,1,1)$ are in
one-to-one correspondence with the set of divisor classes in 
$\Theta-\bigl\{0,X_{1},X_{2}\bigr\}$.
\endproclaim

By symmetry there are completely analogous results for the other three corners.
i.e. the cases $\boldz=(1,0,1,1)$, $\boldz=(1,1,0,1)$, 
and $\boldz=(1,1,1,0)$ correspond to divisor classes in the sets
$\Theta_{X_{1}}-\bigl\{X_{1},0,X_{1}+X_{2}\bigr\}$,
$\Theta_{X_{2}}-\bigl\{X_{2},0,X_{1}+X_{2}\bigr\}$, and
$\Theta_{X_{1}+X_{2}}-\bigl\{X_{1}+X_{2},X_{1},X_{2}\bigr\}$, respecively.

The remainder of the infinity types require more explicit analysis.
There should be no multi-projective points where exactly two of the
$z_{ij}$ are non-zero, nor should there be any multi-projective points
where all the $z_{ij}$ are zero. Where three of the $z_{ij}$ are zero
and one is non-zero, should correspond to a unique multi-projective
point. By symmetry, it suffices to consider the cases
$\boldz=(1,1,0,0)$, $\boldz=(1,0,0,1)$, $\boldz=(1,0,0,0)$, and
$\boldz=(0,0,0,0)$, where
$\boldz=\bigl(z_{\zz},z_{\zo},z_{\oz},z_{\oo}\bigr)$.

\subhead
The Case $\boldz=(1,0,0,0)$
\endsubhead

This really corresponds to four different cases, which can be pictorially
viewed as
$$\matrix
\labeledsquare{z_{\zz}=1}{z_{\zo}=0}{z_{\oz}=0}{z_{\oo}=0}{\scriptstyle(+X_{1})}{\scriptstyle(+X_{2})}{\scriptstyle(+X_{1})}{\scriptstyle(+X_{2})} &
\labeledsquare{z_{\zz}=0}{z_{\zo}=1}{z_{\oz}=0}{z_{\oo}=1}{\scriptstyle(+X_{1})}{\scriptstyle(+X_{2})}{\scriptstyle(+X_{1})}{\scriptstyle(+X_{2})} \cr
\vphantom{\Bigg\vert} \cr
\labeledsquare{z_{\zz}=0}{z_{\zo}=0}{z_{\oz}=1}{z_{\oo}=0}{\scriptstyle(+X_{1})}{\scriptstyle(+X_{2})}{\scriptstyle(+X_{1})}{\scriptstyle(+X_{2})} &
\labeledsquare{z_{\zz}=0}{z_{\zo}=0}{z_{\oz}=0}{z_{\oo}=1}{\scriptstyle(+X_{1})}{\scriptstyle(+X_{2})}{\scriptstyle(+X_{1})}{\scriptstyle(+X_{2})} \cr
\endmatrix
$$
and now by symmetry it is sufficient to consider just the case $z_{\oz}=z_{\zo}=z_{\zz}=0$
and $z_{\oo}=1$, and as before each of the arrows refers
to a simple bihomogeneous case.

Intutively, for this to be possible, this would correspond to points
$Z\in\JacC-\Theta$ such that $Z+X_{1}\in\Theta$ and
$Z+X_{2}\in\Theta$ and
$Z+X_{1}+X_{2}\in\JacC-\Theta$. Now if
$Z\in\JacC-\Theta$ such that $Z+X_{1}\in\Theta$ then $Z$ is 
represented by a reduced divisor of the form
$Q_1+P_{\rtone}-2P_\infty$ for some affine point $Q_1\in C$. Similarly, if
$Z\in\JacC-\Theta$ such that $Z+X_{2}\in\Theta$ then $Z$ is 
represented by a reduced divisor of the form
$Q_2+P_{\rttwo}-2P_\infty$ for some affine point $Q_2\in C$. This would
mean that $Q_1=P_{\rttwo}$ and $Q_2=P_{\rtone}$, so that $Z$ is represented 
by the unique reduced divisor $P_{\rtone}+P_{\rttwo}-2P_\infty$. But then
$Z=X_{1}+X_{2}$ (a 2-division point on $\JacC$), and 
$Z+X_{1}+X_{2}=0$ so $Z+X_{1}+X_{2}\in\Theta$. Thus there
should be a unique point with this infinity type.

The analysis here starts out similar to the previous case.
The conditions $z_{\zz}=1$ and $z_{\zo}=0$ imply that
adding $X_{1}$ gives 
$$\eqalign{
0&= u_{\zz,0}+{\rtone}\,u_{\zz,1}+{\rtone}^2\cr
0&=-v_{\zo,0}+(u_{\zz,1}+{\rtone})\,v_{\zo,1} \cr
0&= u_{\zo,0}\cr
0&= u_{\zo,1}\cr
}$$
while the conditions $z_{\zz}=1$ and $z_{\oz}=0$ imply that 
adding $X_{2}$ gives
$$\eqalign{
0&= u_{\zz,0}+{\rttwo}\,u_{\zz,1}+{\rttwo}^2\cr
0&=-v_{\oz,0}+(u_{\zz,1}+{\rttwo})\,v_{\oz,1} \cr
0&= u_{\oz,0}\cr
0&= u_{\oz,1}\cr
}$$
and since ${\rtone}\ne{\rttwo}$, this allows the values of
$u_{\zz,0}$ and $u_{\zz,1}$ to be determined as
$$\eqalign{
{\rtone}\,{\rttwo}&=u_{\zz,0} \cr
-{\rtone}-{\rttwo}&=u_{\zz,1} \cr
}$$
and therefore
$$\eqalign{
0&=-v_{\zo,0}-{\rttwo}\,v_{\zo,1} \cr
0&=-v_{\oz,0}-{\rtone}\,v_{\oz,1} \cr
}$$
and now note that this forces
$$\eqalign{
0&\ne-v_{\zo,0}-{\rtone}\,v_{\zo,1} \cr
0&\ne-v_{\oz,0}-{\rttwo}\,v_{\oz,1} \cr
}$$
since ${\rtone}\ne{\rttwo}$ since
otherwise either all the $\bolds_{\zo}$ variables 
would all be 0 or all the $\boldS_{\oz}$ variables 
would all be 0.

Now note that 
the conditions $z_{\oo}=0$ and $z_{\oz}=0$ imply that
adding $X_{1}$ gives 
$$\eqalign{
0&=u_{\oo,1} \cr
0&=u_{\oz,0}\,u_{\oo,0} \cr
0&=(v_{\oo,0}+{\rtone}\,v_{\oo,1})\,(v_{\oz,0}+{\rtone}\,v_{\oz,1}) \cr
}$$
and similarly, the conditions $z_{\oo}=0$ and $z_{\zo}=0$ imply that
adding $X_{2}$ gives 
$$\eqalign{
0&=u_{\oo,1} \cr
0&=u_{\zo,0}\,u_{\oo,0} \cr
0&=(v_{\oo,0}+{\rttwo}\,v_{\oo,1})\,(v_{\zo,0}+{\rttwo}\,v_{\zo,1}) \cr
}$$
and since 
$$\eqalign{
0&\ne v_{\zo,0}+{\rtone}\,v_{\zo,1} \cr
0&\ne v_{\oz,0}+{\rttwo}\,v_{\oz,1} \cr
}$$
this implies that 
$$\eqalign{
0&=v_{\oo,0}+{\rtone}\,v_{\oo,1} \cr
0&=v_{\oo,0}+{\rttwo}\,v_{\oo,1} \cr.
}$$
Again, using the fact that ${\rtone}\ne{\rttwo}$ gives
$$\eqalign{
0&=z_{\oo} \cr
0&=u_{\oo,1} \cr
0&=v_{\oo,0} \cr
0&=v_{\oo,1} \cr
}$$
and therefore $u_{\oo,0}\ne 0$ giving rise to single projective
point in the $\boldS_{\oo}$ coordinates, as desired. Note that this
is consistent with $0=u_{\zo,0}\,u_{\oo,0}=u_{\oz,0}\,u_{\oo,0}$ since
it already has been established that $0=u_{\zo,0}=u_{\oz,0}$.

\proclaim{Lemma}
There is exactly one solution to the system of equations $B_4$ with $\boldz=(1,0,0,0)$.
This corresponds to the point $0\in\Pic$.
\endproclaim

By symmetry, the cases $\boldz=(0,1,0,0)$, $\boldz=(0,0,1,0)$, 
and $\boldz=(0,0,0,1)$ also correspond to single points, none
of which have been included in any of the other cases. These
correspond to the points $X_{1}$, $X_{2}$, and $X_{1}+X_{2}$,
respectively.

This exhausts all the divisor classes in $\Pic$. However what has to be checked is whether
there are any ``extra'' points on the quadri-homogeneous variety.

\subhead
The Case $\boldz=(1,0,1,0)$
\endsubhead

This really corresponds to two different cases which can be pictorially viewed as
$$
\labeledsquare{z_{\zz}=1}{z_{\zo}=0}{z_{\oz}=0}{z_{\oo}=1}{\scriptstyle(+X_{1})}{\scriptstyle(+X_{2})}{\scriptstyle(+X_{1})}{\scriptstyle(+X_{2})} 
\labeledsquare{z_{\zz}=0}{z_{\zo}=1}{z_{\oz}=1}{z_{\oo}=0}{\scriptstyle(+X_{1})}{\scriptstyle(+X_{2})}{\scriptstyle(+X_{1})}{\scriptstyle(+X_{2})} 
$$
and now by symmetry it is sufficient to consider just the case $z_{\zz}=z_{\zz}=1$
and $z_{\oz}=z_{\zo}=0$, and as before each of the arrows refers
to a simple bihomogeneous case.

Intutively, if this were possible, there would be some point
$Z\in\JacC-\Theta$ such that $Z+X_{1}\in\Theta$ and
$Z+X_{2}\in\JacC-\Theta$ and
$Z+X_{1}+X_{2}\in\Theta$. This is impossible because if
$Z\in\JacC-\Theta$ and $Z+X_{1}\in\Theta$, then $Z$ is represented
by a divisor of the form $Q+P_{\rtone}-2P_\infty$ for some affine
point $Q\in C$.  Similarly if $Z+X_{2}\in\JacC-\Theta$ and
$Z+X_{1}+X_{2}\in\Theta$, then $Z+X_{2}$ is represented
by a divisor of the form $Q^\prime+P_{\rtone}-2P_\infty$ for some affine
point $Q^\prime\in C$. However $Z+X_{2}$ is also represented by
$Q+P_{\rtone}+P_{\rttwo}-3P_\infty$. Thus
$$Q+P_{\rtone}+P_{\rttwo}-3P_\infty \sim Q^\prime+P_{\rtone}-2P_\infty$$
so
$$Q+P_{\rttwo}-2P_\infty \sim Q^\prime-P_\infty$$
giving a point in $\JacC$ that is both on the $\Theta$-divisor and not on the
$\Theta$-divisor, which of course is impossible. This is to be reflected in the
equations themselves.

The relations between the $\boldS_{\zz}$ and $\boldS_{\zo}$ variables give
$$\eqalign{
0&= u_{\zz,0}+{\rtone}\,u_{\zz,1}+{\rtone}^2\cr 
0&= v_{\zz,0}+{\rtone}\,v_{\zz,1} \cr
0&=u_{\zo,0} \cr
0&=u_{\zo,1} \cr
0&=-v_{\zo,0} + (u_{\zz,1}+{\rtone})\,v_{\zo,1} \cr
}$$
since $z_{\zz}=1$ and $z_{\zo}=0$,
and similarly, the relations between the $\boldS_{\oz}$ and $\boldS_{\oo}$ variables give
$$\eqalign{
0&= u_{\oz,0}+{\rtone}\,u_{\oz,1}+{\rtone}^2\cr 
0&= v_{\oz,0}+{\rtone}\,v_{\oz,1} \cr
0&=u_{\oo,0} \cr
0&=u_{\oo,1} \cr
0&=-v_{\oo,0} + (u_{\oz,1}+{\rtone})\,v_{\oo,1} \cr
}$$
since $z_{\oz}=1$ and $z_{\oo}=0$,
while the relations between the $\boldS_{\zo}$ and $\boldS_{\oo}$ variables give
$$\eqalign{
0&=u_{\zo,0}\,u_{\oo,0} \cr
0&=-(v_{\zo,0}+{\rttwo}v_{\zo,1})\,u_{\oo,0} + (v_{\oo,0}+{\rttwo}v_{\oo,1})\,u_{\zo,0} \cr
0&=(v_{\oo,0}+{\rttwo}v_{\oo,1})(v_{\zo,0}+{\rttwo}v_{\zo,1}) \cr
}$$
since $z_{\zo}=1$ and $z_{\oo}=0$.
Either $v_{\oo,0}+{\rttwo}v_{\oo,1}=0$ or $v_{\zo,0}+{\rttwo}v_{\zo,1}=0$.
Now if $0=v_{\oo,0}+{\rttwo}v_{\oo,1}$, then combined with 
$v_{\oo,0} = (u_{\oz,1}+{\rtone})\,v_{\oo,1}$ gives
$0=(u_{\oz,1}+{\rtone}+{\rttwo})\,v_{\oo,1}$
and noting that $v_{\oo,1}$ can't be 0 (if $v_{\oo,1}=0$, then $v_{\oo,0}=0$, which
would mean that all the $\boldS_{\oo}$ variables were 0). Therefore
$u_{\oz,1}=-{\rtone}-{\rttwo}$. However, $u_{\oz,0}+{\rtone}\,u_{\oz,1}+{\rtone}^2=0$
which then implies that $u_{\oz,0}={\rtone}\,{\rttwo}$. Thus the $U$-polynomial
in the $\boldS_{\oz}$ variables is 
$$U(x)=x^2+u_{\oz,1}\,x+u_{\oz,0}
      =x^2-({\rtone}+{\rttwo})+{\rtone}\,{\rttwo}
      =(x-{\rtone})\,(x-{\rttwo})
$$
and since $\rtone$ and $\rttwo$ are both roots of $f(x)$, this means that
the $V$ polynomial in the $\boldS_{\oz}$ variables, 
$V(x)=v_{\oz,1}\,x+v_{\oz,0}$, must be zero when evaluated
at $x=\rtone$ and at $x=\rttwo$. Therefore $v_{\oz,1}=v_{\oz,0}=0$ 
since $\rtone\ne\rttwo$.

Now consider the glue equations between the $\bolds_{\zz}$ variables and
the $\boldS_{\oz}$ variables. In particular, there is the homogeneous equation
$$\eqalign{0
&=G_1(\bolds_{\zz}M_{\rttwo},\boldS_{\oz}M_{\rttwo},\bolda_{\rttwo}) \cr
&=(u_{\zz,0}+{\rttwo}u_{\zz,1}+{\rttwo}^2z_{\zz})(u_{\oz,0}+{\rttwo}u_{\oz,1}+{\rttwo}^2z_{\oz}) 
- a_1z_{\zz}z_{\oz} \cr}$$
which for $z_{\zz}=z_{\oz}=1$ gives
$$a_1=(u_{\zz,0}+{\rttwo}u_{\zz,1}+{\rttwo}^2)\,(u_{\oz,0}+{\rttwo}u_{\oz,1}+{\rttwo}^2).$$
However,
$U(x)=x^2+u_{\oz,1}\,x+u_{\oz,0}=(x-{\rtone})\,(x-{\rttwo})$,
so ${\rttwo}^2+{\rttwo}u_{\oz,1}+u_{\oz,0}=0$, and therefore $a_1=0$, which is impossible.

In a completely similar fashion, $v_{\zo,0}+{\rttwo}v_{\zo,1}=0$ implies $a_1=0$, as well. 
Therefore the case $\boldz=(z_{\zz},z_{\zo},z_{\oz},z_{\oo})=(1,0,1,0)$ has
no solutions (i.e. there are no such points with this infinity type).

\proclaim{Lemma}
There are no solutions to the system of equations $B_4$ with $\boldz=(1,0,1,0)$.
\endproclaim

By symmetry, the cases $\boldz=(1,1,0,0)$, $\boldz=(0,1,0,1)$, and $\boldz=(0,0,1,1)$
also have no solutions.

\subhead
The Case $\boldz=(1,0,0,1)$
\endsubhead

This really corresponds to four different cases, which can be pictorially
viewed as
$$\matrix
\labeledsquare{z_{\zz}=1}{z_{\zo}=1}{z_{\oz}=0}{z_{\oo}=0}{\scriptstyle(+X_{1})}{\scriptstyle(+X_{2})}{\scriptstyle(+X_{1})}{\scriptstyle(+X_{2})} &
\labeledsquare{z_{\zz}=0}{z_{\zo}=1}{z_{\oz}=0}{z_{\oo}=1}{\scriptstyle(+X_{1})}{\scriptstyle(+X_{2})}{\scriptstyle(+X_{1})}{\scriptstyle(+X_{2})} \cr
\vphantom{\Bigg\vert} \cr
\labeledsquare{z_{\zz}=0}{z_{\zo}=0}{z_{\oz}=1}{z_{\oo}=1}{\scriptstyle(+X_{1})}{\scriptstyle(+X_{2})}{\scriptstyle(+X_{1})}{\scriptstyle(+X_{2})} &
\labeledsquare{z_{\zz}=1}{z_{\zo}=0}{z_{\oz}=1}{z_{\oo}=0}{\scriptstyle(+X_{1})}{\scriptstyle(+X_{2})}{\scriptstyle(+X_{1})}{\scriptstyle(+X_{2})} \cr
\endmatrix
$$
and now by symmetry it is sufficient to consider just the case $z_{\zo}=z_{\zz}=0$
and $z_{\oz}=z_{\oo}=1$, and as before each of the arrows refers
to a simple bihomogeneous case.

Intutively, if this were possible, there would be some point
$Z\in\JacC-\Theta$ such that $Z+X_{1}\in\Theta$ and
$Z+X_{2}\in\Theta$ and
$Z+X_{1}+X_{2}\in\JacC-\Theta$. This is impossible because if
$Z\in\JacC-\Theta$ such that $Z+X_{1}\in\Theta$ then $Z$ is 
represented by a reduced divisor of the form
$P_1+P_{\rtone}-2P_\infty$ for some affine point $P_1\in C$. Similarly, if
$Z\in\JacC-\Theta$ such that $Z+X_{2}\in\Theta$ then $Z$ is 
represented by a reduced divisor of the form
$P_2+P_{\rttwo}-2P_\infty$ for some affine point $P_2\in C$. This would
mean that $P_1=P_{\rttwo}$ and $P_2=P_{\rtone}$, so that $Z$ is represented 
by the unique reduced divisor $P_{\rtone}+P_{\rttwo}-2P_\infty$. But then
$Z=X_{1}+X_{2}$ (a 2-division point on $\JacC$), and 
$Z+X_{1}+X_{2}=0$ so $Z+X_{1}+X_{2}\in\Theta$.

The conditions $z_{\zz}=1$ and $z_{\zo}=0$ imply that
adding $X_{1}$ gives 
$$\eqalign{
0&= u_{\zz,0}+{\rtone}\,u_{\zz,1}+{\rtone}^2\cr
0&=-v_{\zo,0}+(u_{\zz,1}+{\rtone})\,v_{\zo,1} \cr
0&= u_{\zo,0}\cr
0&= u_{\zo,1}\cr
}$$
while the conditions $z_{\zz}=1$ and $z_{\oz}=0$ imply that 
adding $X_{2}$ gives
$$\eqalign{
0&= u_{\zz,0}+{\rttwo}\,u_{\zz,1}+{\rttwo}^2\cr
0&=-v_{\oz,0}+(u_{\zz,1}+{\rttwo})\,v_{\oz,1} \cr
0&= u_{\oz,0}\cr
0&= u_{\oz,1}\cr
}$$
and since ${\rtone}\ne{\rttwo}$, this allows the values of
$u_{\zz,0}$ and $u_{\zz,1}$ to be determined as
$$\eqalign{
{\rtone}\,{\rttwo}&=u_{\zz,0} \cr
-{\rtone}-{\rttwo}&=u_{\zz,1} \cr
}$$
and therefore
$$\eqalign{
0&=-v_{\zo,0}-{\rttwo}\,v_{\zo,1} \cr
0&=-v_{\oz,0}-{\rtone}\,v_{\oz,1}. \cr
}$$

In a completely similar fashion,
the conditions $z_{\oo}=1$ and $z_{\oz}=0$ imply that
adding $X_{1}$ gives 
$$\eqalign{
0&= u_{\oo,0}+{\rtone}\,u_{\oo,1}+{\rtone}^2\cr
0&=-v_{\oz,0}+(u_{\oo,1}+{\rtone})\,v_{\oz,1} \cr
0&= u_{\oz,0}\cr
0&= u_{\oz,1}\cr
}$$
while the conditions $z_{\oo}=1$ and $z_{\zo}=0$ imply that 
adding $X_{2}$ gives
$$\eqalign{
0&= u_{\oo,0}+{\rttwo}\,u_{\oo,1}+{\rttwo}^2\cr
0&=-v_{\zo,0}+(u_{\oo,1}+{\rttwo})\,v_{\zo,1} \cr
0&= u_{\zo,0}\cr
0&= u_{\zo,1}\cr
}$$
and since ${\rtone}\ne{\rttwo}$, this allows the values of
$u_{\oo,0}$ and $u_{\oo,1}$ to be determined as
$$\eqalign{
{\rtone}\,{\rttwo}&=u_{\oo,0} \cr
-{\rtone}-{\rttwo}&=u_{\oo,1} \cr
}$$
and therefore
$$\eqalign{
0&=-v_{\oz,0}-{\rttwo}\,v_{\oz,1} \cr
0&=-v_{\zo,0}-{\rtone}\,v_{\zo,1}. \cr
}$$

However, the conditions
$$\eqalign{
0&=-v_{\zo,0}-{\rttwo}\,v_{\zo,1} \cr
0&=-v_{\zo,0}-{\rtone}\,v_{\zo,1} \cr
}$$
now have the unique solution
$$\eqalign{
0&=v_{\zo,0} \cr
0&=v_{\zo,1} \cr
}$$
since ${\rtone}\ne{\rttwo}$. This now forces
all the $\bolds_{\zo}$ variables to be 0, which
is impossible. A similar argument would force
all the $\bolds_{\zo}$ variables to be 0, as well
(though this is not needed). Therefore there
are no solutions of $B_4$ with $\boldz=(z_{\zz},z_{\zo},z_{\oz},z_{\oo})=(1,0,0,1)$, 
i.e. there are no such points with this infinity type.

\proclaim{Lemma}
There are no solutions to the system of equations $B_4$ with $\boldz=(1,0,0,1)$.
\endproclaim

By symmetry, the case $\boldz=(0,1,1,0)$ also has no solutions.

\subhead
The Case $\boldz=(0,0,0,0)$
\endsubhead

This final case refers to the picture
$$\labeledsquare{z_{\zz}=0}{z_{\zo}=0}{z_{\oz}=0}{z_{\oo}=0}{\scriptstyle(+X_{1})}{\scriptstyle(+X_{2})}{\scriptstyle(+X_{1})}{\scriptstyle(+X_{2})}$$
where all the points at all the corners refer to points at infinity in the projective completions
of the affine components.

First note that $z_{\zz}=z_{\zo}=z_{\oz}=z_{\oo}=0$ also imply
$u_{\zz,1}=u_{\zo,1}=u_{\oz,1}=u_{\oo,1}=0$ and the four sets
of glue equations imply only the following twelve equations
$$\eqalign{
0&=u_{\zz,0}\,u_{\zo,0} \cr
0&=(v_{\zo,0}+{\rtone}v_{\zo,1})\,u_{\zz,0} + (v_{\zz,0}+{\rtone}v_{\zz,1})\,u_{\zo,0} \cr
0&=(v_{\zz,0}+{\rtone}v_{\zz,1})(v_{\zo,0}+{\rtone}v_{\zo,1}) \cr
0&=u_{\zz,0}\,u_{\oz,0} \cr
0&=(v_{\oz,0}+{\rttwo}v_{\oz,1})\,u_{\zz,0} + (v_{\zz,0}+{\rttwo}v_{\zz,1})\,u_{\oz,0} \cr
0&=(v_{\zz,0}+{\rttwo}v_{\zz,1})(v_{\oz,0}+{\rttwo}v_{\oz,1}) \cr
0&=u_{\oz,0}\,u_{\oo,0} \cr
0&=(v_{\oo,0}+{\rtone}v_{\oo,1})\,u_{\oz,0} + (v_{\oz,0}+{\rtone}v_{\oz,1})\,u_{\oo,0} \cr
0&=(v_{\oz,0}+{\rtone}v_{\oz,1})(v_{\oo,0}+{\rtone}v_{\oo,1}) \cr
0&=u_{\zo,0}\,u_{\oo,0} \cr
0&=(v_{\oo,0}+{\rttwo}v_{\oo,1})\,u_{\zo,0} + (v_{\zo,0}+{\rttwo}v_{\zo,1})\,u_{\oo,0} \cr
0&=(v_{\zo,0}+{\rttwo}v_{\zo,1})(v_{\oo,0}+{\rttwo}v_{\oo,1}) \cr
}$$
which are all that follow directly from the above set of glue equations.
In fact all these equations can be derived from just the bilinear part
of the glue.

Now suppose that one of the $u_{ij,0}\ne0$, say $u_{\zz,0}\ne0$ (with the
other cases all being essentially the same). Then it follows that
$u_{\zo,0}=u_{\oz,0}=0$, and from this it follows that 
$(v_{\zo,0}+{\rtone}v_{\zo,1})\,u_{\zz,0}=(v_{\oz,0}+{\rttwo}v_{\oz,1})\,u_{\zz,0}=0$
and therefore $v_{\zo,0}+{\rtone}v_{\zo,1}=v_{\oz,0}+{\rttwo}v_{\oz,1}=0$.
In order that all the $\bolds_{\zo}$ variables not all be zero, it is necessary
that $v_{\zo,0}+{\rttwo}v_{\zo,1} \ne0$ since $\rtone\ne\rttwo$. Similarly
in order that all the $\boldS_{\oz}$ variables not all be zero, it is necessary
that $v_{\oz,0}+{\rtone}v_{\oz,1} \ne0$. This implies (in two different ways)
that $u_{\oo,0}=0$, and also that $v_{\oo,0}+{\rtone}v_{\oo,1}=0$ and
$v_{\oo,0}+{\rttwo}v_{\oo,1}=0$, which then imply that $v_{\oo,0}=v_{\oo,1}=0$
since $\rtone\ne\rttwo$. This gives that all the $\boldS_{\oo}$ variables
are zero, which is impossible. Therefore $u_{\zz,0}=0$.  Similarly
$u_{\zo,0}=u_{\oz,0}=u_{\oo,0}=0$, as well.

With all the $u_{ij,0}=0$, the twelve equations reduce to the following four:
$$\eqalign{
0&=(v_{\zz,0}+{\rtone}v_{\zz,1})(v_{\zo,0}+{\rtone}v_{\zo,1}) \cr
0&=(v_{\zz,0}+{\rttwo}v_{\zz,1})(v_{\oz,0}+{\rttwo}v_{\oz,1}) \cr
0&=(v_{\oz,0}+{\rtone}v_{\oz,1})(v_{\oo,0}+{\rtone}v_{\oo,1}) \cr
0&=(v_{\zo,0}+{\rttwo}v_{\zo,1})(v_{\oo,0}+{\rttwo}v_{\oo,1}). \cr
}$$

Now suppose that $v_{\zz,0}+{\rtone}v_{\zz,1}=0$. In order that the
$\bolds_{\zz}$ variables not all be zero it follows that $v_{\zz,0}+{\rttwo}v_{\zz,1}\ne0$
and therefore $v_{\oz,0}+{\rttwo}v_{\oz,1}=0$. In order that the
$\boldS_{\oz}$ variables not all be zero it follows that $v_{\oz,0}+{\rtone}v_{\oz,1}\ne0$
and therefore $v_{\oo,0}+{\rtone}v_{\oo,1}=0$. In order that the
$\boldS_{\oo}$ variables not all be zero it follows that $v_{\oo,0}+{\rttwo}v_{\oo,1}\ne0$
and therefore $v_{\zo,0}+{\rttwo}v_{\zo,1}=0$. In order that the
$\boldS_{\zo}$ variables not all be zero it follows that $v_{\zo,0}+{\rtone}v_{\zo,1}\ne0$.
This then gives rise to a point on the multi-projective variety.

Alternatively, suppose that $v_{\zz,0}+{\rtone}v_{\zz,1}\ne0$. Then 
$v_{\zo,0}+{\rtone}v_{\zo,1}=0$ and in order that $\bolds_{\zo}$ variables 
not all be zero it follows that $v_{\zo,0}+{\rttwo}v_{\zo,1}\ne0$. Then
$v_{\oo,0}+{\rttwo}v_{\oo,1}=0$ and in order that $\boldS_{\oo}$ variables 
not all be zero it follows that $v_{\oo,0}+{\rtone}v_{\oo,1}\ne0$. Then
$v_{\oz,0}+{\rtone}v_{\oz,1}=0$ and in order that $\boldS_{\oz}$ variables 
not all be zero it follows that $v_{\oz,0}+{\rttwo}v_{\oz,1}\ne0$. Then
$v_{\zz,0}+{\rttwo}v_{\zz,1}=0$. This then gives rise to a second point
on the multi-projective variety.

To summarize, if $z_{\zz}=z_{\zo}=z_{\oz}=z_{\oo}=0$, then the glue equations
still allow the existence of two points. This is a big problem, since
these two points should not exist.

\proclaim{Proposition}
There are exactly two solutions to the system of equations $B_4$ with $\boldz=(0,0,0,0)$.
\endproclaim

Thus working only with two affine Weierstrass points does not lead to
a set of equations for all of $\Jac$. To proceed further additional
equations are needed to eliminate these two points.

\head
Introducing Another Weierstrass Point
\endhead

One approach to resolving this whole problem is to introduce a third
rational Weierstrass point into the analysis.  Each rational Weierstrass point
$P_{\rho^{(i)}}=\bigl(\rho^{(i)},0\bigr)\in\Caff$ is associated with a
point $X_{\rho^{(i)}}=\Cl(P_{\rho^{(i)}}-P_\infty)\in\JacC$ of order 2,
which is also on $\Theta$. Now instead of four copies of the affine
closure of $\JacC-\Theta$, there are eight such copies, and instead of
four sets of inhomogeneous and homogeneous variables, there are now
eight such sets. It is convenient to index them as
$\bolds_{ijk}=\bigl(u_{ijk,0},u_{ijk,1},v_{ijk,0},v_{ijk,1}\bigr)$ and
$\boldS_{ijk}=\bigl(u_{ijk,0},u_{ijk,1},v_{ijk,0},v_{ijk,1},z_{ijk}\bigr)$,
respectively for $ijk\in\I$ with index set
$\I=\{\zzz,\zzo,\zoz,\zoo,\ozz,\ozo,\ooz,\ooo\}$.  There are also
eight copies of the affine variety $A(\bolds)$ that represents
$\JacC-\Theta$ denoted $A(\bolds_{ijk})$ for $ijk\in\I$. The product of
the eight projective closures $\prod_{ijk\in\I}\bar A(\boldS_{ijk})$
is a multi-projective variety in $(\P^4)^8$ and under the Segre
imbedding is a projective variety in $\P^{5^8-1}=\P^{390624}$, a very
large dimensional ambient space.

These projective varieties $\bar A(\boldS_{ijk})$ may be thought of as
corners of a cube, with adjacent corners being pairs of indices of
Hamming distance 1.  The edges of the cube are determined by glue
relations and represent the addition of a point of order 2 in $\JacC$
determined by $X_{\rho^{(1)}}$, $X_{\rho^{(2)}}$, or $X_{\rho^{(3)}}$.
Theese three points determine a subgroup of $\JacC$ of order eight.
Elements of this group will be written as $X_{ijk}$ for $ijk\in\I$.
Here $i$, $j$, and $k$ can take the values $0$ or $1$
$$X_{ijk} 
= k\cdot X_{\rho^{(1)}} + j\cdot X_{\rho^{(2)}} + i\cdot X_{\rho^{(3)}}$$
with a slight abuse of notation.
Thus $X_{ijk}\in\Theta$ if and only if $ijk$ has Hamming weight 0 or 1.
Note that the divisor of the function $y$ on $C$ is given by
$$\div(y)
=P_{\rho^{(1)}}+P_{\rho^{(2)}}+P_{\rho^{(3)}}+P_{\rho^{(4)}}+P_{\rho^{(5)}}
-5\cdot P_\infty$$
so $X_{\rho^{(1)}}+X_{\rho^{(2)}}+X_{\rho^{(3)}}=X_{\rho^{(4)}}+X_{\rho^{(5)}}$.
Each corner of the cube also has the associated affine equations 
$A(\bolds_{ijk})$ which should be thought of as representing 
$\JacC-\Theta_{X_{ijk}}$ where $\Theta_{X_{ijk}}$ is the translate of the
$\Theta$-divisor by the point $X_{ijk}$ of order 2. Each edge of the
cube is associated with the set of bihomogeneous equations between the
variable sets on the two corners at the ends of the edge. These twelve
sets of bihomogeneous equations corresponding to the twelve edges of the
cube define a subvariety of the multi-projective variety
$\prod_{ijk\in\I}\bar A(\boldS_{ijk})$. It is this subvariety, which
is itself a multi-projective variety, that will be identified as $\Jac(C)$.

An infinity type is an association of a value of $1$ or $0$ to each
$z_{ijk}$ on the eight corners of the cube, which determines whether
or not a point in $\JacC$ is an affine point on
$A(\bolds_{ijk})=\JacC-\Theta_{X_{ijk}}$.  Every point in $\JacC$ has some
one or more indices where $z_{ijk}=1$. The question is whether every point in 
multi-projective space defined by the product of the closures $\bar
A(\bolds_{ijk})$ along with the glue equations corresponds to a point
in $\JacC$. In principle, there are $2^{8}=256$ such infinity types,
however, there are a lot of symmetries that reduce the number of cases
that need to be considered. Pictorially, an infinity type is an assignment
of $0$ or $1$ to each corner of the cube
$$\cube{z_{000}}{z_{001}}{z_{010}}{z_{011}}{z_{100}}{z_{101}}{z_{110}}{z_{111}}$$
and it will be important to determine exactly which infinity types are
possible.

It is important to realize that if $ijk$ and $i^\prime j^\prime
k^\prime$ are two indices that are Hamming distance 1 apart from each
other (i.e.  they are connected by an edge of glue relations) and if
$z_{ijk}=0$ while $z_{i^\prime j^\prime k^\prime}=1$, then the
previous analysis (between the hatted and the unhatted variables with
$z=1$ and $\zhat=0$) shows that the $u_0$, $u_1$, $v_0$, and $v_1$
values for $i^\prime j^\prime k^\prime$ completely determine the
$u_0$, $u_1$, $v_0$, and $v_1$ values for $ijk$.

It is useful to look at parts of the cube, where some of the vertices
and their corresponding edges are considered.  All the discussions
that can be restricted to faces of the cube are derivable from the
quadri-projective analysis.  There are three cases of high interest that
are parts of the cube that live in three dimensions.  

The first case of interest is an analysis of a corner and its three
edges where the corner will have $z=1$ but the
vertices at the opposite ends of each of the three edges will have
$z=0$.  Thus the case is $z_{000}=1$ and $z_{001}=z_{010}=z_{100}=0$, i.e.
$$
\cubecorner{z_{000}}{z_{010}}{z_{001}}{z_{100}}
\qquad=\qquad
\cubecorner{1}{0}{0}{0}
$$
with the goal of showing that there are no points with this partial 
infinity type. Consider the three faces that intersect at this corner.
Each of these faces has at least two $0$'s and at least one $1$. However,
the square analysis shows that it is impossible to have two $1$'s and 
two $0$'s, therefore the corners need to be $1$ and the infinity type
of the cube then must look like
$$
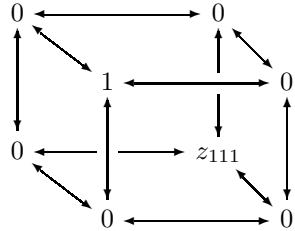
\cube{1}{0}{0}{0}{0}{0}{0}{z_{111}}
$$
with $z_{111}$ undetermined. Now each face is a square of type
$(z_{\zz},z_{\oz},z_{zo},z_{\oo})=(1,0,0,0)$, which was analyzed
above. Having three $0$'s in a square turns out to be very restricting.
If the order 2 points being added in the glue along the edges of the square
are $X_{\rho^{(1)}}$ and $X_{\rho^{(2)}}$, then $U_{\zz}(x)
=\bigl(x-\rho^{(1)}\bigr)\,\bigl(x-\rho^{(2)}\bigr)$ and $V_{\zz}(x)=0$.
Similarly, if the two points are $X_{\rho^{(1)}}$ and $X_{\rho^{(3)}}$
then $U_{\zz}(x)=\bigl(x-\rho^{(1)}\bigr)\,\bigl(x-\rho^{(3)}\bigr)$ 
and if the two points are $X_{\rho^{(2)}}$ and $X_{\rho^{(3)}}$
then $U_{\zz}(x)=\bigl(x-\rho^{(2)}\bigr)\,\bigl(x-\rho^{(3)}\bigr)$.
However, $\rho^{(1)}$, $\rho^{(2)}$, and $\rho^{(3)}$ are all distinct
values. This shows that this infinity type cannot exist.
This same analysis will apply to any corner of the cube to show that
the analogous partial infinity type cannot exist.

The second case of interest is an analysis of a corner and its three
edges where the corner will have $z=0$ and the
vertices at the opposite ends of each of the three edges will also have
$z=0$.  Thus the case is $z_{000}=z_{001}=z_{010}=z_{100}=0$, i.e.
$$
\cubecorner{z_{000}}{z_{010}}{z_{001}}{z_{100}}
\qquad=\qquad
\cubecorner{0}{0}{0}{0}
$$
and the goal now is to show that for the three faces that come
together at this corner, the opposite vertices must all have
$z=1$, i.e. that $z_{110}\ne0$, $z_{101}\ne0$, and $z_{011}\ne0$.
Consider first what happens if one of the opposite vertices has a
1 and one of the other opposite vertices has a 0. Thus the picture 
(omitting the $z_{111}$ corner) is
$$
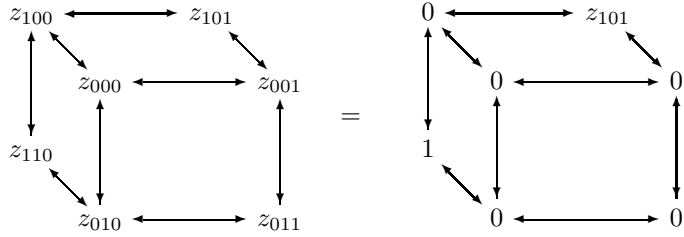
\cubeminuscorner{z_{000}}{z_{001}}{z_{010}}{z_{011}}{z_{100}}{z_{101}}{z_{110}}
\qquad=\qquad
\cubeminuscorner{0}{0}{0}{0}{0}{z_{101}}{1}
$$
where the assignments $z_{011}=0$ and $z_{110}=1$ have been made (and
the other possibilities are equivalent by symmetry, so the same
analysis will apply). Now consider the all zero face
$\{000,001,010,011\}$.  The quadri-projective (two-dimensional)
analysis showed that if this were to occur then the $u_0$ variable at
the $000$ vertex would have to be $u_{000,0}=0$.  Now the face with a
single 1 (the $\{000,100,010,110\}$ face also was analyzed in the
quadri-projective case, and there it was shown that the $u_0$ variable
at the $000$ vertex would have to be $u_{000,0}\ne0$. This is obviously
a contradiction, so the assignments $z_{011}=0$ and $z_{110}=1$ are
invalid, and this analysis applies regardless of whether, $z_{101}=0$
or $z_{101}=1$.

Next consider what happens if all three of the opposite vertices have
$z=0$, i.e. 
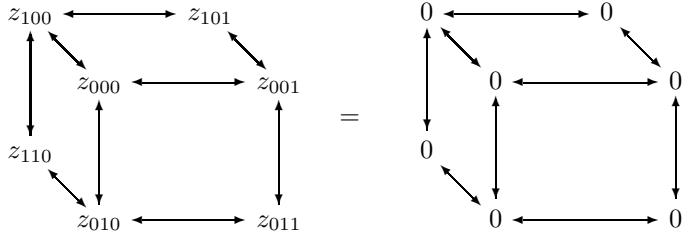
$$\hskip-25pt
\cubeminuscorner{z_{000}}{z_{001}}{z_{010}}{z_{011}}{z_{100}}{z_{101}}{z_{110}}
\qquad=\qquad
\cubeminuscorner{0}{0}{0}{0}{0}{0}{0}
$$
again omitting the $z_{111}$ corner from the picture. Now the prior
quadri-projective (two-dimensional) analysis of
$(z_{\zz},z_{\zo},z_{\oz},z_{\oo})=(0,0,0,0)$ applies to each of the
three faces. In particular, it must be the case that at the $\zzz$
vertex, $u_{\zzz,0}=0$. Furthermore, the top face
$\{000,001,010,011\}$, which involves adding $X_{\rho^{(1)}}$ and
$X_{\rho^{(2)}}$, requires that either
$v_{\zzz,0}+\rho^{(1)}v_{\zzz,1}=0$ or
$v_{\zzz,0}+\rho^{(2)}v_{\zzz,1}=0$, but not both (or else all the
$\boldS_{\zzz}$ will be $0$. Similarly, the $\{000,001,100,101\}$
face, which involves adding $X_{\rho^{(1)}}$ and $X_{\rho^{(3)}}$,
requires that either $v_{\zzz,0}+\rho^{(1)}v_{\zzz,1}=0$ or
$v_{\zzz,0}+\rho^{(3)}v_{\zzz,1}=0$, but not both, and the
$\{000,010,100,110\}$ face, which involves adding $X_{\rho^{(2)}}$ and
$X_{\rho^{(3)}}$, requires that either
$v_{\zzz,0}+\rho^{(2)}v_{\zzz,1}=0$ or
$v_{\zzz,0}+\rho^{(3)}v_{\zzz,1}=0$, but not both. These conditions
cannot all be simultaneously satisfied.  If $A$, $B$, and $C$ are
logical propositions, and if $((A \OR B) \AND (A \OR C) \AND (B \OR
C))$ is to be true, then at least two out the three logical
propositions $\{A,B,C\}$ must hold. This then shows that
$z_{011}=z_{101}=z_{110}=0$ cannot hold. Therefore it must be the
case that 
$$
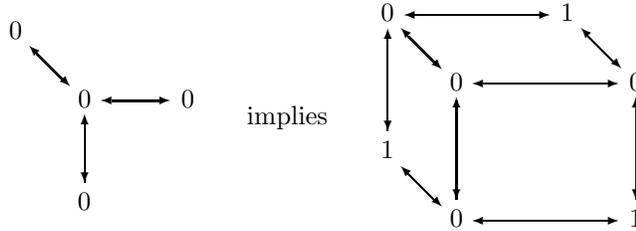
\cubecorner{0}{0}{0}{0}
\qquad\hbox{implies}\qquad
\cubeminuscorner{0}{0}{0}{1}{0}{1}{1}
$$
and to finish off the $z_{111}$ corner, note that that there cannot be
two $0$'s and two $1$'s on any face, and therefore $z_{111}=1$,
i.e. the full infinity type is
$$
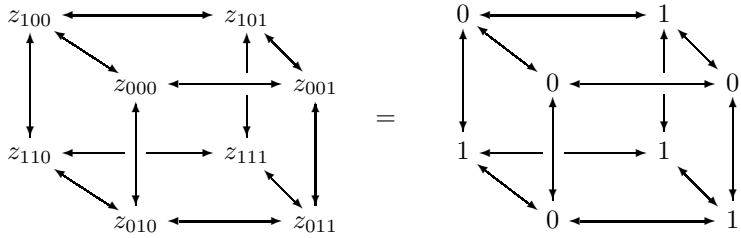
\cube{z_{000}}{z_{001}}{z_{010}}{z_{011}}{z_{100}}{z_{101}}{z_{110}}{z_{111}}
\qquad=\qquad
\cube{0}{0}{0}{1}{0}{1}{1}{1}
$$
In this case it is possible to be more specific and state precisely
what is at each corner.

Note that another way of saying this is that any ball of Hamming
radius 2 must contain at least one $z=1$ value. In particular, this
approach shows that in three dimensions (i.e by using three affine
Weierstrasss points), no face can have $z=0$ at all its corners.  This
result could not be shown by staying only in two dimensions by using
only two affine Weierstrass points.

It is interesting that this infinity type corresponds to a unique point on
$\JacC$. This follows from the two-dimensional case for infinity type
$\{z_{\zz},z_{\zo},z_{\oz},z_{\zo}\}=\{0,0,0,1\}$ where now three of 
the four $1$'s on the cube are determined to correspond to the points 
$X_{\rho^{(1)}}+X_{\rho^{(2)}}$, $X_{\rho^{(1)}}+X_{\rho^{(3)}}$,
$X_{\rho^{(2)}}+X_{\rho^{(3)}}$, all of which are in $\JacC-\Theta$.
The fourth point at the $111$ vertex corresponds to 
$X_{\rho^{(1)}}+X_{\rho^{(2)}}+X_{\rho^{(3)}}
=X_{\rho^{(4)}}+X_{\rho^{(5)}}$ which is also in $\JacC-\Theta$.

The third case of interest is an analysis of a corner and its three
edges where the corner will have $z=0$ but the
vertices at the opposite ends of each of the three edges will have
$z=1$.  Thus the case is $z_{000}=0$ and $z_{001}=z_{010}=z_{100}=1$, i.e.
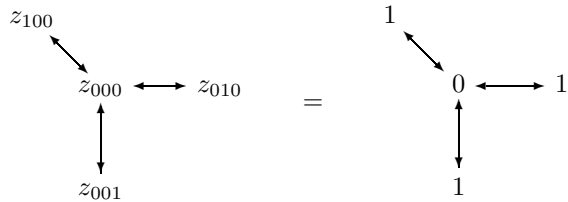
$$
\cubecorner{z_{000}}{z_{010}}{z_{001}}{z_{100}} 
\qquad=\qquad
\cubecorner{0}{1}{1}{1}
$$
and now note that each of the three faces that contain this corner has two $1$'s
and one $0$. This means that the far corners on these three faces must all be $1$'s
with the final $111$ vertex to be determined, i.e.
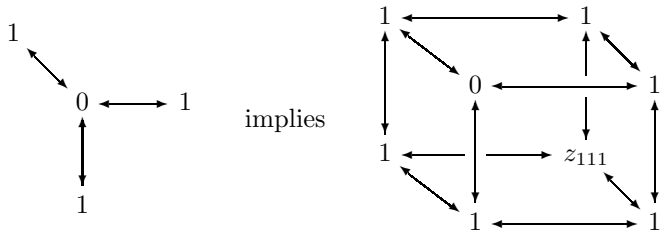
$$
\cubecorner{0}{1}{1}{1}
\qquad\hbox{implies}\qquad
\cube{0}{1}{1}{1}{1}{1}{1}{z_{111}}
$$
with $z_{111}$ to be determined.
There are thus two possible infinity types
$$
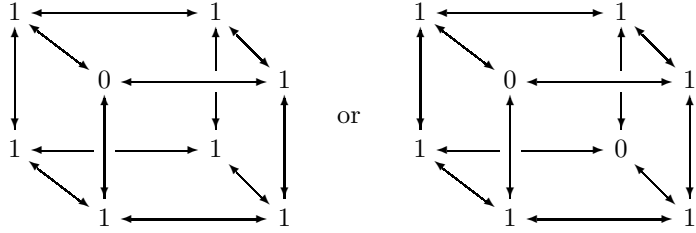
\cube{0}{1}{1}{1}{1}{1}{1}{1}
\qquad\hbox{or}\qquad
\cube{0}{1}{1}{1}{1}{1}{1}{0}
$$
which can be analyzed more completely. 

In the first case, where there is only one $z=0$ component, these are
points on $\Pic$ of the form $P-P_{\infty}$ where $P$ is an affine
point that is not a Weierstrass point on the curve $C$.  In the second
case, the situation corresponds to $P$ being either $P_4$ or $P_5$,
i.e. one of the two Weierstrass points on $C$ that have not been part
of the analysis up to now. Since $P_1+P_2+P_3+P_4+P_5-5\cdot
P_{\infty}$ is the divisors of the function $y$ on $C$, if $P=P_4$ or
$P=P_5$, the points at the $z_{\zzz}$ corner will be on $\Theta$ and
the points at the $z_{\ooo}$ corner will also be on $\Theta$ (actually
on $\Theta_{X_1+X_2+X_3}$).  The way to see this is as follows.
Suppose the reduced divisor at $z_{\zzz}=0$ is $P-P_{\infty}$ and the
reduced divisor at $z_{\ooo}=0$ is $P^\prime-P_{\infty}$, adding
$X_1+X_2+X_3$ moves from one to the other, i.e. the divisor class
represented by the non-reduced divisor $P^\prime+P_1+P_2+P_3-4\cdot
P_{\infty}$ is the same as the divisor class represented by
$P-P_{\infty}$, so these divisors are equivalent.  This means that the
divisors $P^\prime+\iota(P)-2\cdot P_{\infty}$ and $P_1+P_2+P_3-3\cdot
P_{\infty}$ are equivalent. However $P_1+P_2+P_3-3\cdot P_{\infty}$ is
equivalent to $P_4+P_5-2\cdot P_{\infty}$, and therefore the divisors
$P^\prime+\iota(P)-2\cdot P_{\infty}$ and $P_4+P_5-2\cdot P_{\infty}$
are equivalent, but these are both reduced, unless it happens that
$P^\prime=P$. However if $P^\prime=P$ then the reduction of the
divisor $P^\prime+\iota(P)-2\cdot P_{\infty}$ is $0$, which certainly
not the same as $P_4+P_5-2\cdot P_{\infty}$.  Therefore there is
actual equality as divisors, i.e.  $P^\prime+\iota(P)-2\cdot
P_{\infty} = P_4+P_5-2\cdot P_{\infty}$, which means that either
$P=P_4$ and $P^\prime=P_5$ or else $P=P_5$ and $P^\prime=P_4$.

\proclaim{Proposition}
For the infinity type where two antipodal corners (say $(\zzz)$ and $(\ooo)$) 
of the cube have $z_{\zzz}=z_{\ooo}=0$ and all other
cube corners have $z_{ijk}=1$, there are exactly two multi-projective points.
These correspond to the points $X_4$ and $X_5$ of order 2 in $\Pic$.
\endproclaim

If this analysis had included all the Weierstrass points, these two
points would not look special. However, that would have corresponded
to a 4-dimensional cube (a tesseract), rather than a 3-dimensional
cube, in which the 16 corners would be in one-to-one correspondence
with the 16 points of order 2 on $\JacC$. Fortunately, there is no
need to go to these lengths to construct the Jacobian.

It is interesting to consider in three dimensions what happens if
$ijk$ and $i^\prime j^\prime k^\prime$ are two indices that are
Hamming distance 1 apart with both $z_{ijk}=0$ and $z_{i^\prime
j^\prime k^\prime}=0$. Now such an edge is common to two different
faces of the cube. Each face of the cube is then subject to the
analysis for the quadri-projective case.  It was shown there that on
any square face, it is impossible to have an infinity type consisting
of two $1$'s and two $0$'s. Now suppose that $z_{000}=z_{100}=0$. The
two faces to consider are $\{000,100,101,001\}$ and
$\{000,100,110,010\}$, and there has to be at least one more zero on
each of these two faces. Up to symmetry, there are two possibilities:
$$
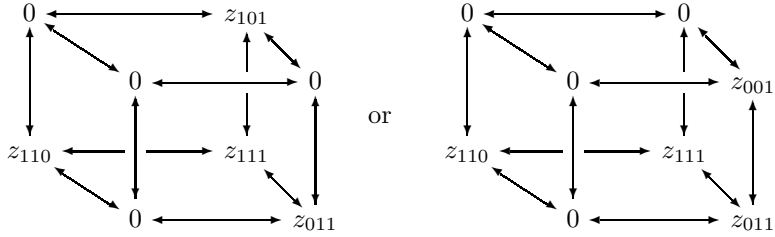
\cube{0}{0}{0}{z_{011}}{0}{z_{101}}{z_{110}}{z_{111}}
\qquad\hbox{or}\qquad
\cube{0}{z_{001}}{0}{z_{011}}{0}{0}{z_{110}}{z_{111}}
$$
i.e. either $z_{010}=z_{101}=0$ or $z_{001}=z_{010}=0$. The first
possibility has already been studied above and leads to a unique completion
and a unique point on $\JacC$. Now consider the second
possibility.  There are now two faces $\{000,001,011,010\}$ and
$\{100,101,110,111\}$ that have two zeros, so now there has to be
another zero on each of these two faces. The analysis above shows
that it is necessary that $z_{001}=z_{110}=1$ or else the $000$ corner
would have it and all its adjacent vertices and at least one vertex
of Hamming distance 2 from $000$ with $z=1$. Thus there is the partial
infinity type
$$
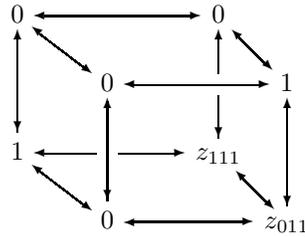
\cube{0}{1}{0}{z_{011}}{0}{0}{1}{z_{111}}
$$
and now note that the infinity type completes to
$$
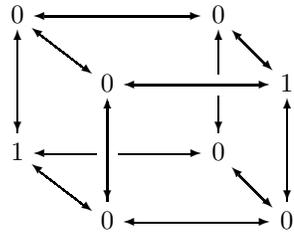
\cube{0}{1}{0}{0}{0}{0}{1}{0}
$$
since no face can have two $0$'s and two $1$'s. However, this partial
infinity type is also 
impossible. To see this, Note that the
$\{000,001,011,010\}$ face requires the $u_0$ variable at the $010$
vertex to be $u_{\zoz,1}\ne0$, while the $\{000,010,110,100\}$ face
requires the $u_0$ variable at the $010$ vertex to be $u_{\zoz,1}=0$. 

Therefore, if some edge of the cube, say the $\{000,100\}$ edge,
has partial infinity type $(0,0)$, i.e. $z_{000}=z_{100}=0$, there
are only two possible valid ways to complete the infinity type, either
$$
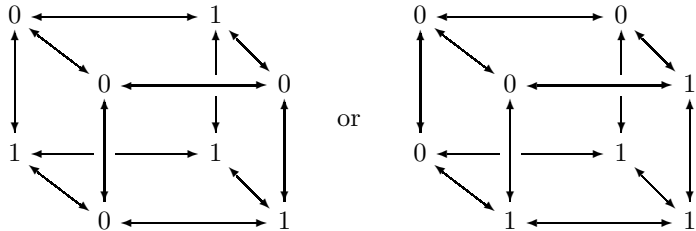
\cube{0}{0}{0}{1}{0}{1}{1}{1}
\qquad\hbox{or}\qquad
\cube{0}{1}{1}{1}{0}{0}{0}{1}
$$
which are equivalent (up to symmetry), each of which determines
a single point in $\JacC$.

\proclaim{Proposition}
If some edge of the cube, say the $\{000,100\}$ edge,
has partial infinity type $(0,0)$, i.e. $z_{000}=z_{100}=0$
then one of the two following two cases holds, each of which 
corresponds to a single point in $\JacC$.
\vskip 1pt (i) 
$({z_{000}},{z_{001}},{z_{010}},{z_{011}},{z_{100}},{z_{101}},{z_{110}},{z_{111}})
=({0},{0},{0},{1},{0},{1},{1},{1})$ corresponding to the divisor class 0, represented
by the reduced divisor 0;
\vskip 1pt (ii)
$({z_{000}},{z_{001}},{z_{010}},{z_{011}},{z_{100}},{z_{101}},{z_{110}},{z_{111}})
=({0},{1},{1},{1},{0},{0},{0},{1})$ corresponding to the divisor class $X_1$, represented
by the reduced divisor $P_1-P_{\infty}$.
\endproclaim

Basically, if there are two corners of the cube with $z=0$, the
situation is very constrained, and there are only finitely many
possibilities for which points on $\JacC$ they correspond to. More
importantly, there are no anomalous points.  All the divisor classes
in $\Pic$ correspond to unique points on the multi-projective variety
and every multi-projective point on the variety corresponds to a
unique reduced divisor. The $\Theta$-divisor and all its translates by
$[2]$-division points that are some combination of $X_1$, $X_2$, and
$X_3$ correspond to the set of points at infinity on each of the
corners of the cube in the projective closures that are at each of the
affine varieties on the corners.

\subhead
Summary of the Octo-Homogeneous Case
\endsubhead

This completes the construction of $\JacC$. It requires that only
three of the affine Weierstrass points be rational over the base field.
The Jacobian is constructed as a multi-projective variety sitting inside
of the multi-projective space $(\P^4)^8$, and becomes a projective variety
under the Segre map.  The points on $\JacC$ are in one-to-one correspondence
with divisor classes in $\Pic$ that are rational over the base field.
It is useful to summarize this as follows. 

Let $K$ be a field of characteristic other than 2.
Let $f(x)$ be a monic quintic polynomial over $K$ having no repeated roots and having at least three
roots in $K$, i.e. write
$$\eqalign{f(x) 
&= x^5+a_4\,x^4+a_3\,x^3+a_2\,x^2+a_1\,x+a_0 \cr
&=\bigl(x-\rho^{(1)}\bigr)\,\bigl(x-\rho^{(2)}\bigr)\,\bigl(x-\rho^{(3)}\bigr) 
 \,\bigl(x-\rho^{(4)}\bigr)\,\bigl(x-\rho^{(5)}\bigr) \cr
}$$
with $a_0,a_1,a_2,a_2,a_3,a_4,\rho^{(1)},\rho^{(2)},\rho^{(3)}\in K$ and 
write 
$$\bolda=\bigl(a_0,a_1,a_2,a_2,a_3,a_4\bigr).$$ 
For any $\rho\in K$ define
the matrices
$$M(\rho)=
\pmatrix 1       & 0     & 0    & 0      & 0 \cr
        \rho     & 1     & 0    & 0      & 0 \cr
        0        & 0     & 1    & 0      & 0 \cr
        0        & 0     & \rho & 1      & 0 \cr
        {\rho}^2 & 2\,\rho & 0      & 0      & 1 \cr
\endpmatrix
\qquad\hbox{and}\qquad
A(\rho) = 
\pmatrix
1      & 0         & 0         & 0       & 0 \cr
\rho   & 1         & 0         & 0       & 0 \cr
\rho^2 & 2\,\rho   & 1         & 0       & 0 \cr
\rho^3 & 3\,\rho^2 & 3\,\rho   & 1       & 0 \cr
\rho^4 & 4\,\rho^3 & 6\,\rho^2 & 4\,\rho & 1 \cr
\endpmatrix
$$
and the vector
$$
\boldb(\rho)
=
\pmatrix
\rho^5 & 5\,\rho^4 & 10\,\rho^3 & 10\,\rho^2 & 5\,\rho \cr
\endpmatrix
.$$

For a vector of variables $\boldS=\bigl(u_0,u_1,v_0,v_1,z\bigr)$ let
$$\eqalign{
E_0(\boldS,\bolda) 
&= a_0\,z^4 - a_2\,u_0\,z^3 + a_4\,u_0^2\,z^2 + a_3\,u_0\,u_1\,z^2 
 - 2\,u_0^2\,u_1\,z - a_4\,u_0\,u_1^2\,z + u_0\,u_1^3 \cr&\quad
 - v_0^2\,z^2 + u_0\,v_1^2\,z \cr
E_1(\boldS,\bolda) 
&= a_1\,z^4 - a_3\,u_0\,z^3 + u_0^2\,z^2 - a_2\,u_1\,z^3 + 2\,a_4\,u_0\,u_1\,z^3 
 + a_3\,u_1^2\,z^2 - 3\,u_0\,u_1^2\,z \cr&\quad - a_4\,u_1^3\,z + u_1^4 
 - 2\,v_0\,v_1\,z^2 + u_1\,v_1^2\,z^2 \cr
}$$
and for convenience write $\boldE(\boldS,\bolda)=\bigl(E_0(\boldS,\bolda),E_1(\boldS,\bolda)\bigr)$.
For a second vector of variables $\boldhatS=\bigl(\uhat_0,\uhat_1,\vhat_0,\vhat_1,\zhat\bigr)$ let
$$\eqalign{
G_1(\boldS,\boldhatS,\bolda) &= u_0\,\hat u_0 - a_1\,z\,\hat z \cr
G_2(\boldS,\boldhatS,\bolda) &= \hat v_0\,u_0 + v_0\,\hat u_0 \cr
G_3(\boldS,\boldhatS,\bolda) &= v_0\,\hat v_0 + a_1\left(a_4\,z\,\hat z  - u_1\,\hat z - \hat u_1\,z\right) \cr
G_4(\boldS,\boldhatS,\bolda) &= \hat u_0\left(v_1\,\hat z+\hat v_1\,z\right) + \hat v_0\left(u_1\,\hat z - \hat u_1\,z\right) \cr
G_5(\boldS,\boldhatS,\bolda) &= u_0\left(v_1\,\hat z+\hat v_1\,z\right) + v_0\left(\hat u_1\,z - u_1\,\hat z\right) \cr
G_6(\boldS,\boldhatS,\bolda) &= a_1\,\hat u_0\,z^2 - a_1\,a_3\,z^2\,\hat z + a_1\,u_0\,\hat z\,z + a_1\,\hat u_1\,u_1\,z -
2\,u_1\,\hat v_0\,v_0 - 2\,\hat u_0\,v_0\,v_1 \cr
G_7(\boldS,\boldhatS,\bolda) &= a_1\,\hat u_0\,z\,\hat z - a_1\,a_3\,z\,\hat z^2 + a_1\,u_0\,\hat z^2 + a_1\,\hat u_1\,u_1\,\hat z -
2\,\hat u_1\,\hat v_0\,v_0 - 2\,u_0\,\hat v_0\,\hat v_1 \cr
G_8(\boldS,\boldhatS,\bolda) &= -a_1\,\hat v_0\,z^3 + a_3\,u_0\,\hat v_0\,z^2 - u_0^2\,\hat v_0\,z -
2\,a_4\,u_0\,u_1\,\hat v_0\,z + 3\,u_0\,u_1^2\,\hat v_0 
\cr&\quad + a_1\,u_1\,\hat v_1\,z^2 +
a_1\,u_1\,v_1\,z\,\hat z + 2\,\hat v_0\,v_0\,v_1\,z \cr
G_9(\boldS,\boldhatS,\bolda) &= -a_1\,v_0\,\hat z^3 + a_3\,\hat u_0\,v_0\hat z^2 - \hat u_0^2\,v_0\hat z -
2\,a_4\,\hat u_0\,\hat u_1\,v_0\,\hat z + 3\,\hat u_0\,\hat u_1^2\,v_0 
\cr&\quad +a_1\,\hat u_1\,v_1\,\hat z^2 + a_1\,\hat u_1\,\hat v_1\,\hat z\,z + 2\,v_0\,\hat v_0\,\hat v_1\,\hat z \cr
}$$
and for convenience write $\boldG(\boldS,\boldhatS,\bolda)
=\bigl\{G_1(\boldS,\boldhatS,\bolda),\ldots,G_9(\boldS,\boldhatS,\bolda)\bigr\}$.

Let 
$$\I=\Z/2 \times \Z/2 \times \Z/2$$ 
be an index set, and let
$$\eqalign{
\boldd_1&=(1,0,0)\cr
\boldd_2&=(0,1,0)\cr
\boldd_3&=(0,0,1)\cr
}$$
be a basis (mod 2)of $\I$. For each $\boldc=(i,j,k)\in\I$ let
$$\boldS_{\boldc}
  =\bigl(u_{\boldc,0},u_{\boldc,1},v_{\boldc,0},v_{\boldc,1},z_{\boldc}\bigr)$$ 
be a vector of projective variables. The elements of $\I$ correspond to 
the corners of a cube, and each edge of the cube is $(\boldc,\boldc+\boldd_l)$
for some $\boldc\in\I$ and $l\in\{1,2,3\}$ (but note that in this notation
the edges $(\boldc,\boldc+\boldd_l)$ and $(\boldc+\boldd_l,\boldc)$ are the same).

Then for each $\boldc\in\I$ (i.e. at each corner of the cube) 
there are the defining homogeneous equations 
$$\boldE(\boldc; \bolda)
=\bigl\{E_0(\boldS_{\boldc},\bolda),E_1(\boldS_{\boldc},\bolda)\bigr\}$$
and for each $\boldc\in\I$ and each $l\in\{1,2,3\}$ (i.e. at each edge of the cube)
there are the defining bihomogeneous equations
$$\boldG(\boldc,l;\bolda) 
= \left\{G_i(M\bigl(\rho^{(l)}\bigr)\boldS_{\boldc},
         M\bigl(\rho^{(l)}\bigr)\boldS_{\boldc+\boldd_l},
         A(\rho^{(l)})\bolda+\boldb(\rho^{(l)}))\right\}_{i=1}^{i=9}
$$
and note that each edge has been counted twice in this notation, 
so each glue equation appears twice. There are a total 16 homogeneous equations
for the corners and 108 bihomogeneous equations for the edges.

\proclaim{Theorem}
Let $C$ be the genus 2 hyperelliptic curve
$$y^2=f(x)$$ 
with $$f(x)= x^5+a_4\,x^4+a_3\,x^3+a_2\,x^2+a_1\,x+a_0$$
with coefficients in a field $K$ not of characteristic 2, and such that the
$f(x)$ has no multiple roots and at least three roots in $K$. 
Let $J$ be the multi-homogeneous
variety defined by the polynomials $\boldE(\boldc; \bolda)$ as $\boldc$ ranges
over $\I$ and also by the polynomials $\boldG(\boldc,l;\bolda)$ as $\boldc$ ranges
over $\I$ and $l$ ranges over $\{1,2,3\}$. Then $J$ is a multi-homogeneous model
of $\JacC$.
\endproclaim

\enddocument